\documentclass{svmult}

\usepackage[bottom]{footmisc}
\usepackage{amsmath}
\usepackage{amssymb,latexsym}
\usepackage[all]{xy}

\makeindex

\newtheorem{lbremark}{{\bf Remark}}[section]

\def\bee{\begin{equation}}
\def\ee{\end{equation}}
\def\ga{\gamma}
\def\ka{\kappa}
\def\lam{\lambda}
\def\th{\theta}
\def\Ga{\Gamma}
\def\de{\delta}
\def\eps{\epsilon}
\def\De{\Delta}
\def\de{\delta}
\def\om{\omega}
\def\Om{\Omega}

\def\kc{\mathcal{K}}

\def\pc{{\mathcal{P}}}
\def\gc{{\mathcal{G}}}
\def\pa{{P}^{\mathrm{ad}}}
\def\pac{{\mathcal{P}^{\mathrm{ad}}}}
\def\osc{{\mathcal{O}_S}}
\def\oxc{{\mathcal{O}_X}}
\def\la{\longrightarrow}
\def\LG{\mathrm{Lie}\,(G)}
\def\ot{\otimes}
\def\fg{\mathfrak{g}}
\def\fgx{\mathfrak{g}_X}
\def\pra{p_{0}^{\ast}}
\def\prb{p_{1}^{\ast}}
\def\prc{{p_{2}^{\ast}}}
\def\prd{{p_{3}^{\ast}}}

\def\dea{{\Delta^{1}_{X/S}}}
\def\deb{{\Delta^{2}_{X/S}}}
\def\dec{{\Delta^{3}_{X/S}}}

\def\noi{\noindent}
\def\te{\vartheta}
\def\wtg{\widetilde{\gamma}}

\begin{document}

\title*{Differential Geometry of Gerbes\\ and Differential Forms}
\author{Lawrence Breen\thanks{Unit\'e Mixte de Recherche   CNRS 7539}}

\institute{Institut Galil\'ee\\
Universit\'e Paris 13\\
99, avenue J. - B. Cl\'ement
\\ 93430 Villetaneuse, France\\
\texttt{breen@math.univ-paris13.fr}}

\maketitle

\numberwithin{equation}{section}

\bigskip

\begin{center}
{\it 
To Murray Gerstenhaber and Jim Stasheff
}
\end{center}

\begin{abstract}
We discuss certain aspects of the combinatorial approach to the
differential geometry of non-abelian gerbes due to W. Messing and the
author \cite{dgg}, and give a more direct
derivation of the associated cocycle equations. This
 leads us to  a more
  restrictive definition  than in \cite{dgg} of the  corresponding 
 coboundary relations. We also show that the
 diagrammatic proofs of certain local curving and curvature equations may
 be replaced by computations with differential forms.
\end{abstract}

\section{Introduction}
\setcounter{footnote}{0}

It is a  classical fact\footnote{at least in a differential geometric setting,
see \cite{kn}, but the same construction can be carried out within the
context of algebraic geometry.} that to a principal $G$-bundle $P$  on a scheme $X$,
endowed with a connection $\epsilon$, is associated a   $\LG$-valued
2-form $\kappa$ on $P$, the curvature of the connection, 
satisfying a certain $G$-equivariance condition. While $\kappa$ does not
in general
descend to a 2-form on $X$, the equivariance
condition may  be viewed   as a descent
condition for $\kappa$ from a 2-form on $P$  to a 2-form on $X$, but
now 
 with values in the Lie algebra of the gauge group $P^{\mathrm{ad}}$ of $P$.
 The connection on $P$ also induces a connection $\mu$ on the group
$\pa$, and the  2-form $\kappa$  satisfies the Bianchi equation, 
an equation  which may  be
expressed in global terms as 
\begin{equation}
\label{dmu}
 \D\kappa + [\mu,\kappa] = 0 
\end{equation}
(\cite{dgg} proposition 1.7, \cite{cdf} theorem 3.7). 
Choosing a local trivialization of the bundle $P$, 
 on an open cover 
$\mathcal{U}:= \coprod_{i \in I} U_i$ of $X$, the connection $\epsilon$ is
described by a family of  $\LG$-valued  connection 1-forms 
$\om_i$ defined on  the open sets
$U_i$, and the associated curvature $\kappa$  corresponds to a family 
of  $\LG$-valued 2-forms $\kappa_i$ defined, according to the 
so-called structural equation of Elie Cartan, by the formula\footnote{The
  canonical divided power $1/2[\om,\, \om]$ of the 2-form $[\om,\,\om]$
is also  denoted $\om \wedge \om$ or $[\om]^{(2)}$.}

\bee
\label{elie}
\kappa_i = \D\om_i + \dfrac{1}{2}\,[\om_i,\, \om_i]
\ee
 Equation \eqref{dmu} then  reduces to 
the classical Bianchi identity
\bee
\label{bianchiclas}
\D\kappa_i + [\om_i, \kappa_i]  = 0\; .
\ee

\bigskip

  J. -L. Brylinski  introduced in  \cite{bry} the notions of 
connection $\epsilon$ and  curving $K$ on  an abelian  $G$-gerbe $\pc$
 on a space
$X$ (where $G$ was the multiplicative group $G_m$, or  rather in his
framework the group  $U(1)$), and  showed that  to such connective data
$(\epsilon,\,K)$ is
associated a closed
$G_m$-valued 3-form $\om$ on $X$, the 3-curvature. More recently,
W. Messing and the author extended these concepts in \cite{dgg}
 from abelian to general, not necessarily 
 abelian, gerbes $\pc$ on a scheme $X$. The coefficients of
such a gerbe  no longer constitute a sheaf of groups as in the
principal
 bundle
situation, but rather a monoidal stack $\gc$ on $X$, as is
to be expected in that categorified setting. In particular, when  the  gerbe is
associated to a given non-abelian group $G$ (so that we refer to it as a
$G$-gerbe), the  corresponding  coefficient stack $\gc$ is the monoidal stack
associated to the prestack determined  by the crossed module $G \la
\mathrm{Aut}(G)$, where $\mathrm{Aut}(G)$ is the sheaf of local
automorphisms of $G$. It may also  be
described more invariantly as the monoidal stack of $G$-bitorsors on
$X$. Once more, to the gerbe $\pc$ is associated its gauge stack, a
 twisted form $\pac:= \mathcal{E}q(\pc,\,\pc)$ of the given monoidal
 stack $\gc$, and the
 connection on $\pc$ induces a connection $\mu$  on $\pac$.
By  analogy with the principal bundle case,  the
corresponding 3-curvature $\Om$,
viewed as a global 3-form on $X$, now  takes its
values in the arrows of  the stack $\pac$.

\bigskip

 There now arises   a new, and at first sight somewhat surprising
 feature, but which is simply another facet of the
categorification context in which we are operating.  The 3-form
 $\Om$ is accompanied by an auxiliary 2-form
$\kappa$ with values in the objects of the gauge  stack $\pac$, which we called in \cite{dgg} the {\it fake curvature} of
the given connective  structure  $(\epsilon,\, K)$. A first relation between
the forms $\Om$ and $\kappa$ comes from the  very definition \cite{dgg} 
(4.1.20), (4.1.22) of $\Om$, and may be stated as in
\cite{dgg}(4.3.8)
 as
 the categorical equation
\begin{equation}
\label{eqom1}
t\Om + \D\kappa + [\mu,\, \kappa] = 0
\end{equation}
where $t$ stands for ``target'' of a 1-arrow with source the
identity object $I$ in the stack of $\mathrm{Lie}(\pac)$-valued 3-forms on
$X$. On the other hand, 
the 3-form $\Om$ is no
longer closed, even in the $\mu$-twisted sense described for principal
bundles by
\eqref{dmu}. It satisfies instead the following more complicated analogue 
  \cite{dgg} (4.1.33)  of the Bianchi identity \eqref{dmu}:
 \begin{equation}
\label{eqom2}
\D \Om + [\mu,\, \Om] + [\mathcal{K},\, \kappa] = 0\;.
\end{equation}
 While the  first two terms in this equation are similar to those of \eqref{dmu},  the categorification term  $\mathcal{K}$ is an arrow in the stack of 2-forms with values
 in the  monoidal
 stack $\mathcal{E}q(\pac,\, \pac)$ induced by the curving $K$. 
The  pairing of $\kc$ with $\kappa$
 is induced by the evaluation of the natural transformation $\kc$
 between functors from $\pac$ to itself on the object $\kappa$ 
of $\pac$.

\bigskip

The price to be paid for the compact form in which the  global curvature 
 equations 
\eqref{eqom1} and \eqref{eqom2} have been  stated is their rather
abstract nature, and it is of interest to describe them in a more
local form in terms of traditional group-valued differential forms,
just as  
was done in \eqref{bianchiclas} for equation \eqref{dmu}.
Such a local description  was already obtained in \cite{dgg},  both for the cocycle conditions \eqref{eqom1}
and \eqref{eqom2}, and  for the corresponding coboundary
equations which arise when alternate  local trivializations of the gerbe
have been chosen. However, the determination of those local equations was
 rather indirect, as it required  a third description of a
gerbe, which we have called the semi-local description \cite{ima} \S 4,
 and which has also 
appeared elsewhere in a various situations \cite{ul1}, 
\cite{murray}, \cite{hitch}.

\bigskip

 The present text may be viewed as a companion piece to the author's
 \cite{ima}. 
 Its main  purpose is  to provide a more transparent  construction
 than in
\cite{dgg}
 of the
cocycle  conditions and related equations associated to
 a gerbe with curving data  summarized 
in \cite{dgg}
 theorem 6.4.
  We restrict our attention, as in
 \cite{ima}, to gerbes which are connected rather than   locally
 connected, as these determine    \v{C}ech cohomology classes. A
 cocyclic description in the general
 case requires 
hypercovers  and could be 
 dealt with along the lines discussed in
   \cite{2-gerbe}, but would
not shed any additional light on the phenomena being investigated here.
Our  main results are to be found in sections 4 and 5, while section 3
reviews for the reader's convenience some aspects of \cite{dgg} and
\cite{ima}. Section 2 is a review of some of the formulas in the
differential calculus of $\LG$-valued forms, a few of which 
 do not appear to be well-known.

\bigskip

Another aim of the present work is to revisit the quite
complicated  coboundary equations of \cite{dgg} \S 6.2.  The
 coboundary equations   which arise here are 
  simpler, and  more
consistent than
those of \cite{dgg} with a non-abelian \v{C}ech-de Rham interpretation. We refer to  remark \ref{remark} for a specific comparison
between the two notions. In order to make this comparison easier, we
have chosen the orientations of our arrows consistenly with
\cite{dgg}. This accounts  for example
 for the strange choice of orientation of the
arrow $B_i$ in diagram \eqref{def:bi}, or for the change of sign
\eqref{gatilde} for the arrow $\ga_{ij}$.  
 
\bigskip

A final purpose of this text is to explain how the diagrammatic proofs of
some of the local results of \cite{dgg} can be replaced by more classical
computations involving $\LG$-valued  differential forms. For this
reason,
 we
have given two separate computations for  certain  equations, one
diagrammatic and the other classical. We do not assert that one
 of the two methods of proof is always preferable, though
 one might contend that diagrams
provide  a better understanding of the situation than the
corresponding manipulation of differential forms. As the level of
categorification increases, so will the dimension of
the diagrams to be considered, and it may not be 
 realistic to expect  to tread
 along the diagrammatic path much  beyond the hypercube proof 
\cite{dgg} (4.1.33)
of the
higher Bianchi equation \eqref{eqom2}. The generality and algebraicity
of the formalism of differential forms  must then come into its own. In addition, 
 it is our hope 
 that   the present  approach, which  extends to the gerbe
context  the
traditional methods of differential geometry,
   will provide    an  accessible point of  entry into this topic.
A number of other
authors have recently described  certain  aspects of  the differential geometry of gerbes in
terms
 of differential forms,  particularly  \cite{jurco}, \cite{laurent-st-xu}, and
 \cite{schth}, \cite{baez-schr}. 

\bigskip

 I wish  to thank  Bernard Julia and  Camille Laurent-Gengoux
 for enlightening  discussions on related topics. The impetus for the
 present work was provided by my collaboration with Wiliam Messing on
 our joint papers \cite{cdf} and  \cite{dgg}. It is a pleasure to thank him
 here for our instructive and wide-ranging discussions over all these years.

\bigskip

\section{Group-valued differential forms}
\label{section1}

\subsection{} 

Let $X$ be an $S$-scheme. We assume from now on for simplicity
that 
 that the primes 2 and
3 are invertible in the ring of functions of $S$ (for example $S =
\mathrm{Spec}(k)$ where $k$ is a field of characteristic $\neq 2, 3$).
 A relative differential $n$-form on an $S$-scheme $X$,
with values in a sheaf of  $\osc$-Lie algebras $\fg$  is defined
 as a global section of
the sheaf $\fg \ot_{\osc} \Om^n_{X/S}$ on $X$. When $X/S$ is smooth,
\bee
\label{comp:gforms}
\fg \ot_{\osc} \Om^n_{X/S} \simeq \mathrm{Hom}_{\oxc}(T^n_{X/S},\,
\fgx)
\ee
where $\fgx:= \fg \ot_{\osc}\oxc$ and   $T^n_{X/S}$ is  the $n$-th
exterior power  $\wedge^n T_{X/S}$  of
the relative tangent sheaf $T_{X/S}$, i.e the sheaf of relative
$n$-vector fields on $X$. Such an $n$-form
  is nothing else than an $\oxc$-linear map 
\bee
\label{def:lgform}
T^n_{X/S} \la \fgx\,.
\end{equation}
 
 In view of this definition, such a map
is classically called a $\fg$-valued  differential form.
A more geometric description
of such forms is given in \cite{cdf}, following the
ideas of  A. Kock  in the context of synthetic differential geometry
 \cite{AK:DFG}, \cite{AK:Bianchi}. It is based on the consideration,
 for any positive
 integer $n$, of the scheme $\Delta^n_{X/S}$ of relative infinitesimal
 $n$-simplexes on $X$. For any $S$-scheme $T$, a  $T$-valued point of $\Delta^n_{X/S}$ 
consists of
an $(n+1)$-tuple of $T$-valued points $(x_0,\,\ldots,\, x_n)$ of $X$
which are pairwise close to first order in an appropriate sense
\cite{cdf} (1.4.9). We view $\De^n_X$ as an $X$-scheme via the
projection $p_0$  of such  points to $x_0$.  As $n$ varies, the schemes $\Delta^n_{X/S}$  determine a
 simplicial  $X$-scheme $\Delta^\ast_{X/S}$, whose face and
 degeneracy operations are induced
 by the usual projection and injection
 morphisms $X^n \la X^{n \pm 1}$.

\bigskip

Let $G$ be a flat $S$-group scheme, with $\osc$-Lie algebra $\fg$. 
 A relative $\fg$-valued $n$-form \eqref{def:lgform} on $X/S$ may
then be identified by \cite{cdf} proposition 2.5 
with a morphism of $S$-schemes
\bee
\label{def:nform}
\Delta^n_{X/S} \stackrel{f}{\la} G 
\ee
whose restriction to the degenerate subsimplex $s\De^n_{X/S}$ of
$\De^n_{X/S}$ factors
 through the unit section of $G$. When differential forms are
 expressed in this  combinatorial language, they deserve to be called 
 $G$-valued 
 differential forms, even though they actually coincide with the traditional
 $\fg$-valued differential forms \eqref{comp:gforms},
\eqref{def:lgform}. In the combinatorial context, our notation will
 be multiplicative,
 and additive  when we pass to the traditional language
of  differential forms.

\bigskip
  
 We will now discuss some of the features of these $\fg$-valued
 forms, and  refer to  \cite{cdf} for further discussion.
 First of all, let us recall that the action of the symmetric
 group $S_{n+1}$ on a combinatorial
 differential $n$-form $\om(x_0,\ldots,x_n)$ by permutation of the
 variables is given by
\[
\label{perm}
\om(x_{\sigma(0)},\ldots,x_{\sigma (n)}) = \om(x_0,\ldots,x_n)^{\epsilon(\sigma)}
\]
where $\epsilon(\sigma)$ is the signature of $\sigma$. Also,
 the commutator pairing
\[
[g,\,h] := g\,h\,g^{-1}h^{-1}
\]
 on the group  $G$ determines a
 bracket pairing on $\fg$-valued forms of degree $\geq 1$, defined combinatorially by the rule

\begin{equation}
\label{def:bra}
\xymatrix@C=10pt@R=8pt{
(\fg \ot_\osc \Om^m_{X/S}) \times  (\fg \ot_\osc \Om^n_{X/S})
&\ar[r]^{\:\:\:\:} 
 &&
  (\fg \ot_\osc \Om^{m+n}_{X/S})
}
\end{equation}

\noi which sends $(\om,\,\om')$ to $[\om,\,\om']$, where
\[
 [\om,\,\om'](x_0,\,\ldots, x_{m+n}) :=  [\om (x_0,\, \ldots
 \,,x_m),\, \om'(x_m,\, \ldots,x_{m+n})]\,.
\]

\medskip

\noi This pairing is  defined in classical terms,  by 
\[
\label{bracketdef}
[\om,\,\om']:= [Y,\,Y'] \ot (\eta \wedge \eta')
\]
 for any pair of forms $\om :=Y \ot \eta$ and
$\om':=Y' \ot \eta'$ in  $\fg
\ot_{\osc} \Om^\ast_{X/S}$. 
It  endows $\fg \ot_{\osc} \Om^\ast_{X/S}$ with the structure of a graded
 $\osc$-Lie algebra. In particular, the bracket satisfies
 the graded commutativity
 rule 
\bee
\label{eq:grcom}
[f,\,g] = (-1)^{|f||g|+1}[g,\,f]\;,
\ee 
where $|f|$ is the degree of the form $f$, so that
\[
[f,f] = 0
\] whenever $|f|$ is even. The graded
 Jacobi identity  is expressed (in additive notation) as:
\[
\label{jac}
(-1)^{|f||h|}[f,\,[g,\,h]] + (-1)^{|f||g|}[g,[h,\,f]] +
 (-1)^{|g||h|}[h,\,[f,\,g]] 
 = \, 0\,.
\]
In particular,
\bee
\label{jac1}
[f,\, [f,\,f]] = 0
\ee
and, when $|f| = |g| = 1$, 
\[
\label{jac3}
[f,\, \tfrac{1}{2}\,  [g,\,g]] = [[f,\,g]\,,g]\;.
\] 

\bigskip
Let $\mathrm{Aut}(G)$ be the sheaf of local automorphisms of $G$, whose
group of sections above an $S$-scheme $T$ is the group
$ \mathrm{Aut}_T(G_T)$  of
automorphisms of the $T$-group $G_T:= G \times_S T$. The
definition \eqref{def:nform} of a combinatorial $n$-form still makes
sense when $G$ is replaced by a sheaf of groups  $F$ on $S$, and the
traditional description of such combinatorial $n$-forms as $n$-forms with values in the Lie algebra of $F$
remains valid by \cite{cdf}
proposition 2.3 when $F =
\mathrm{Aut}(G)$. The evaluation map 
\[
\begin{array}{ccc}
\mathrm{Aut}(G) \times G &\la& G\\
(u,\, g) &\mapsto & u(g)
\end{array}\]
induces for all pair of positive integers a bilinear pairing
\bee
\label{def:bra1}
\xymatrix@C=10pt@R=8pt{
(\mathrm{Lie} \,(\mathrm{Aut} (G))\ot_{\osc} \Om^m) 
 \times  (\fg \ot_\osc \Om^n_{X/S})
&
 \ar[r]^{\:\:\:\:} 
&\;\;
  (\fg \ot_\osc \Om^{m+n}_{X/S})
}
\ee
which sends $ (u,\,g)$ to $[u, \, g]$, where
\bee
\label{def:bra-a}
 [u,\,g](x_0,\ldots,x_{m+n}) :=  u(x_0,\ldots,x_m)(g(x_m,\,\ldots,
x_{m+n}))\,\, g(x_m,\,\ldots,
x_{m+n})^{-1}\,.
\ee

\medskip

\noi This pairing is compatible with the  pairings \eqref{def:bra}
associated to the $S$-groups $G$ and $\mathrm{Aut}(G)$ in the following
sense. For any pair of $\fg$-valued forms $g,\,g'$, and  an
  $\mathrm{Aut}(G)$-valued form $u$, 
\begin{equation}
\label{comp:i-bra}
[i(g),\, g'\,] = [g,\,g'\,] \;\;\;\;\; \text{and}\;\; \;\;\;
 i([u,\,g]) = [u,\,
i(g)]
\end{equation}
where $i: G \la \mathrm{Aut}(G)$
 is the inner conjugation map $i(\gamma)(g) :=  \ga\, g \,\ga^{-1}$. 
More generally,  an isomorphism $r: G \la G'$ induces  a morphism $r$
  from $G$-valued combinatorial $n$-forms 
to $G'$-valued combinatorial $n$-forms, compatible with the Lie bracket
operation \eqref{def:bra},  and which  corresponds in
classical terms to the morphism
 $\mathrm{Lie}(r) \ot_{osc} 1: \fg \ot_{\osc}
 \Om^n_{X/S}
\la \fg' \ot  \Om^n_{X/S}$. The functoriality of the
bracket  \eqref{def:bra1} is expressed by the formula
\bee
\label{compbra}
r([u,\, g] = [\,{}^{r}\!u,\, r(g)]
\ee
where ${}^{r\,}\!u := r\,u\,r^{-1}$.
\bigskip

When $u$ is an
 $\mathrm{Aut}(G)$-valued form of degree  $m \geq 1$  and 
 $g$ is a $G$-valued function,  the definition of a pairing 
 \[\begin{array}{ccc} 
(\mathrm{Lie \ Aut}(G) \ot_{\osc} \Om^m_{X/S}) \times  G &\la& \fg
\ot_{\osc} \Om^m_{X/S}\\
(u,\,g) & \mapsto& [u,\, g]
\end{array}
\] is still given by the formula \eqref{def:bra-a}, but
now  with  $n= 0$. This pairing
are no longer linear in $g$, but instead  satisfies the equation
\[
 [u,\, g\,g'] = [u,\, g] + {}^{g\,}\![u,\,g']
\]
where for any $G$-valued form $\om$  and any $G$-valued function 
$g$
the adjoint  left action ${}^{g\,\,}\!\om$ of a function $g$ on a form 
$\om$ is defined
combinatorially by 
\[
({}^{g\,}\!\om)(x_0,\ldots,x_n) := g(x_0) \;\om(x_0,\ldots,x_n)
\; g(x_0)^{-1}\;,
\]
(and  this expression is in fact equal to 
$g(x_i) \;\om(x_0,\ldots,x_n)
\; g(x_i)^{-1}$ for any \linebreak  $0 \leq i \leq n$).  In classical notation
this corresponds,
for $\om = Y \ot \eta\:\, \in \: \mathfrak{g} \ot \Om^n_{X/S}$,
to  the formula
\[
{}^{g\,}\!(Y \ot \eta) =   {}^{g\,\,}\!Y \ot \eta \]
for the adjoint  left action of $g$ on $Y$. The  adjoint right action
${\om}^{\,\ga}$ is defined by
\[{\om}^{\,g}:= {}^{(g^{-1})}\,\om\]
so that 
\[  \om^{g} (x_0,\ldots,x_n) =g(x_0)^{-1}\,
 \om (x_0,\ldots,x_n) \;g(x_0) \;.\]

\bigskip

Similarly, when $g$ is a $G$-valued and $u$ an 
$\mathrm{Aut}(G)$-valued form, a pairing $[g,\,u]$ is defined  by the combinatorial formula
\bee
\label{def:bra1op}
[g,\,u](x_0,\,\ldots,x_{m+n}) := g(x_0,\ldots, x_m) 
\,( u(x_m,\ldots,x_{m+n})(g(x_0,\ldots,x_m)^{-1})) \;. 
\ee 
The pairing \eqref{def:bra1op} satisfies the analogue 
\[
\label{eq:grcom1}
[g,\, u] = (-1)^{|g||u| +1} [u,\,g]
\]
of the graded commutativity rule \eqref{eq:grcom}, so that its
properties may be deduced from those of the pairing $[u,\,g]$. In particular
\[\label{crosshomprop}
[g^{-1},\,u] = - [u,\, g^{-1}] = [u,g]^{g}\;.
\]
\noi We refer to  appendix A of \cite{dgg}  for additional properties of these pairings.

\bigskip

\subsection{}
 
 The de Rham  differential map 
\begin{equation}
\label{deRham}
 \xymatrix@C=40pt{\fg \ot_\osc \Om^n_{X/S} \ar[r]^{d^n_{X/S}} &
\fg \ot_\osc \Om^{n+1}_{X/S}} 
 \end{equation}
is defined combinatorially for $n \geq 2$,
 in Alexander-Spanier fashion,   by 
   \bee
\label{def:dncomb}
\D^n_{X/S}\om(x_0,\,\ldots, x_{n+1}):=
 \prod_{i=0}^{n+1}\om(x_0,\ldots,\widehat{x_i},\ldots
 \om_{n+1})^{(-1)^i}\,.\ee
This definition agrees for $n > 1$ with the  
 classical definition of the
$G$-valued
de Rham differential:
\bee
\label{defdn}
 \D^n_{X/S}\,
 \om := \D_{X/S}\,\, \om
\ee
 where for  $\om = Y \ot \eta\ $ in $\fg
\ot \Om^n_{X/S},\ $

\bee 
\label{defdn1}
\D_{X/S}\,\om:= Y \ot \D\eta \,.\ee

\bigskip

 In particular $\D^n$ is an $\osc$-linear map whenever  $n \geq 2$, 
and  it follows from \eqref{defdn1}  that the composite $d^{n+1} \,
d^n$ is trivial . This also follows from the
combinatorial definition of  $d^n$, since for $n \geq 2$ the factors in
the expression \eqref{def:dncomb} for $d^n\om$ commute with each other.  

\bigskip 

For any section $g$ of $G$, we set
\bee
\label{defdoa}
\D^0_{X/S}(g) := g(x_0)^{-1}g(x_1)\;.
\ee
The map  
\bee
\label{defdo}
\begin{array}{ccc}
G_X &\stackrel{\D^0_{X/S}}{\la} &\fg \ot_{\osc} \Om^1_{X/S}\\ 
g & \mapsto & g^{-1} \D g\,
\end{array}
\end{equation} 
is
a crossed homomorphism, 
for  the  adjoint left  action of $G$ on $\fg$. Observe that the expression 
$g^{-1} \D g$ is consistent with the  combinatorial
definition \eqref{defdoa} of $\D^0_{X/S}(g)$.
 While this  traditional expression  of $\D^0_{X/S}(g)$
 as a product of the  two 
terms  $g^{-1}$ and $\D g$ does  make sense whenever $G$ is a
subgroup scheme of the linear group $GL_{n,S}$, such a decomposition
 is purely conventional for a   general $S$-group scheme $G$.
A companion to $\D^0_{X/S}$ is the differential $\widetilde{\D}^0:  G \la \fg
\ot_{\osc}\Om^1_{X/S}\;$, defined by 
\[
\widetilde{\D}^{\,0}_{X/S}(g)(x_0,\, x_1):= g(x_1) g(x_0)^{-1}\;. \]
The traditional notation for this expression is
$dg \, g^{-1}$.  This notation is consistent with such formulas (in
additive notation)
 as 
\[
\begin{array}{lcccr}
{}^{g\,}\!(g^{-1}dg) = \D g \,g^{-1}& &\text{and}&& -\, (g^{-1}\D g)  =
\D g^{-1}\,g \;.
\end{array}
\]
The differential $\D^1_{X/S}$ is defined
combinatorially by 
\bee
\label{defd1comb}
(\D^1_{X/S}\,\om)(x,y,z):= \om(x,y)\, \om(y,z)\, \om(z,x).
\ee In classical terms, it follows (see \cite{cdf} theorem 3.3)  that 
\bee
\label{defd1}
 \D^1_{X/S} \,\om := \D \om + \frac{1}{2}[\om,\,\om].
 \ee
We will henceforth denote  $\D^n_{X/S}$ simply by $\D^n$ for all $n$.

\bigskip

The    quadratic term    $\frac{1}{2}[\om,\,\om]$ implies that $
\D^1_{X/S}$ is not a linear map, in fact it follows from
\eqref{defd1}, or the elementary combinatorial calculation of
\cite{cdf} lemma 3.2, that
\[
\D^1(\om + \om') = \D^1\om + \D^1\om' +[\om,\, \om']\;.
\]
In particular, 
\[
\label{eq:minom}
\D^1(-\om) = - \D^1(\om) +[\om,\, \om] \;.
\]
 It is immediate, from the combinatorial point of view, that
\bee
\label{d0d1}
\D^1 \D^0(g) =  \D^1(g^{-1}\,\D g) = 0  \ee
 for all $g$ in $G$.
The differential $\D^1$  has a companion, which we will denote by
 $\widetilde{\D}^1$, defined by 
\[\widetilde{\D}^1(\om)(x,\,y,\,z) := \om(z,\, x)\, 
\om(y,\, z)\,\om(x,\, y)\,. \]
A 
 combinatorial computation implies  that
\begin{eqnarray*}
\widetilde{\D}^1 \om &= \D^1 \om - [\om,\, \om]\\
 &\;\;\;= \;\D \om\; - \frac{1}{2} [\om,\, \om]\,,
\end{eqnarray*}
and the analogue
\[
\label{defd0d1a}
\widetilde{\D}^1 (\widetilde{\D}^0(g)) = \widetilde{\D}^1(dg\,g^{-1})
= 0
\]
of \eqref{d0d1} is satisfied.
Finally, it follows from \eqref{defdn} that the  $\D^n$ satisfy  
\[\D^{i+j} [\om,\, \om'] = [\D^i\om, \, \om'] +
(-1)^{i}[\om,\, \D^j\om'] \]
whenever $i,j \geq 2$, and the corresponding formula for the pairing
$[u,\,g]$  \eqref{def:bra-a}
is also valid.

\bigskip

\subsection{}
 We now choose, for any $S$-scheme $X$ and any $S$-group scheme $G$,
  an  $\mathrm{Aut} (G)$-valued 1-form $m$ on $X$.  We  extend
the definition of
 the de Rham differentials \eqref{defdo}, \eqref{defd1comb} and 
\eqref{deRham}
to the twisted  differentials
 \bee
\label{def:mucon-a}
\D^n_{X/S,\,m}:\fg \ot_\osc \Om^n_{X/S} \la
  \fg  \ot_\osc
  \Om^{n+1}_{X/S}
 \end{equation}
(or simply $\D^n_m$) defined combinatorially by the following formulas:
\begin{eqnarray*}
\label{def:d1ma}
\D^1_m \,\om(x_0, x_{1}) & := &
\om(x_0,x_1)\,\, m(x_0,x_1)(\om(x_1,\, x_2)) \,\,m(x_0,x_1)\, m(x_1,x_2)
(\om (x_2,x_0))\\
 &\ =&\om(x_0,x_1)\,\, m(x_0,x_1)(\om(x_1,\, x_2)) \,\om(x_0,\,x_2)^{-1}
\end{eqnarray*}
 \begin{multline*}
\label{def:dnma}
\D^n_m\, \om(x_0,\ldots, x_{n+1}) :=\\
= m(x_0,x_1) (\om(x_1,\ldots x_{n+1}))\,
\prod_{i=1}^{n+1}\om(x_0,\ldots,\widehat{x_i},\ldots,
 x_{n+1})^{(-1)^i} 
\end{multline*}
when
   $n >1$. When the $\mathrm{Aut}(G)$-valued form $m$ is the image $i(\eta)$
under inner conjugation of a $G$-valued form $\eta$, the expression
$\D^n_{i(\eta)} \om$
will  simply be  denoted $\D^n_\eta\, \om$.
The corresponding degree zero  map
$
\D^0_{m}: G \la \fg \ot_{\osc}\Om^1_{X/S}
$ is defined by
\[
\label{def:doma}
\D^0_m (g) := g(x_0)^{-1}\,m(x_0,\, x_1)(g(x_1)),\,
\]
(and $\D^0_m (g)$ 
  will also be  denoted
 $g^{-1}\D_m(g)$, consistenly with \eqref{defdoa}).

\bigskip

 It follows  from elementary combinatorial
computations that  the differentials $\D^n_m$ can  be defined
in classical terms  by
 \bee 
\label{def:dtwist}\D^n_m\om = \D^n\om +[m,\,\om]
\ee for all $n$, so that for any $\fg$-valued 1-form $\eta$,
\bee
\label{addm}
\D^n_{m + i_{\eta}}(\om) = \D^n_m(\om) +[\eta,\,\om]\;.
\ee In particular,
\[
\label{def:dm1}
 \D^1_m(\om) = \D^1 \om +[m,\,\om] = \D\om +\frac{1}{2} [\om,\,\om] +
 [m,\, \om] \,.
\]
While 
the map $\D^n_m$ is linear for $n \geq 2$, 
\bee
\label{d1add}
 \D^1_m(\om + \om') = \D^1_m\om + \D^1_m\om' + [\om,\,\om']
\end{equation}
so that
\bee
\label{d1add-a}
\D^1_m(-\om) = - \D^1_m(\om) - [\om,\,\om].
\end{equation}

\newpage

\noi  Finally, for any section $g$ of $\Gamma$, 
\[
\label{def:dm0}
g^{-1}\D_m g = g^{-1}\D g + [m,\,g]\,.
\]

\bigskip

   The composite  morphism $\D^{n+1}_m\: \D^n_m$ is in
general non-trivial, and the
previous classical definitions of $\D^n_m$  imply that
\bee
\label{dndn1}
\D_m^{n+1} \D_m^n \,\om = [\D^1m,\, \om] 
\ee
whenever $n \geq 2$. For $n= 0$, the corresponding formulas are
\bee
\label{dndn1a}
\begin{array}{ccccc}
\D^1_m\ \D^0_m g = [g^{-1},\,\D^1m] && \text{and} &&
\widetilde{\D}^1_m \ \widetilde{\D}^0_m g = [\D^1m,\, g]
\end{array}
\ee so that, for  $n \neq 1$, we recover the well-known assertion
that the vanishing of   $\D^1m = 0$ implies that 
$\D^{n+1}\D^n = 0$. One verifies that  for any  1-form $\om$
\begin{eqnarray}
\label{falsebianchi}
\D^2_m\, \D^1_m (\om) & =& [\D^1 m, \om]  + [\D^1_m \om,\, \om]
\\
& =& [\D^1m, \om] + [ \D^1 \om,\, \om] + [[m,\,\om],\,\om] 
\;.
\end{eqnarray}
This reduces to   the equation
\[  \D^2_m\, \D^1_m (\om) = [\D^1m,\, \om]\]
of type \eqref{dndn1} whenever
$\D^1_m\om = 0$. For $m = i(\om)$, equation 
\eqref{falsebianchi} is equivalent to 
the classical Bianchi identity  \cite{kn} II Theorem 5.4:
\bee
\label{def:bianchiclas}
\D^2_\om\,\,  \D^1 \om = 0\;.
\ee

\bigskip

We now state the functoriality properties of the differential
\eqref{def:dtwist} $\D^n_m$  for $n \geq 1$. We define the twisted
conjugate ${}^{g\,\ast\,}\!\om$ of a $G$-valued  1-form $\om$ by
 \begin{align}
\label{cocep13clas-b}
  {}^{g\,\ast\,}\!\om := (\pra g)\, \,\om\, \, (\prb 
     g)^{-1}& =  {}^{g\,}\!\om  + g \,
 \D g^{-1}\\
&=\:\, \om + [g,\,\om] + g\, \D g^{-1} \;. \notag
\end{align}
It 
follows from the combinatorial definition \eqref{defd1comb} of $\D^1$
that   
\bee
\label{funct:d1}
 {}^{g\,}\!(\D^1\om)  =  \D^1( {}^{g\,\ast\,}\!\om) \;.
\ee
More generally, for any  $G$-valued form $\om$ of  degree $n \geq 1$,
 and any section
$u$ of $ \mathrm{Aut}(G) $ on $X$, 
\begin{align}
u(\D^n_m(\om)) &= \D^n_{\,({}^{u\,\ast\,}\!m)}\, u(\om) \label{funct:dnm}
\\
&= 
\D^n_{\,({}^{\,u\,}\!m)}\, u(\om) +[u\,\D u^{-1},\, u(\om)]\notag \\
 &= \D^n_m(u(\om)) + [[u,\,m],\,u(\om)] +[u\,\D u^{-1},\,
u(\om)]\;.
\label{funct:dnm1}
\end{align}

\bigskip
 
\section{Gerbes and their connective structures}
\label{section2}

\subsection{}
\label{subsection:2-1}

Let $\pc$ be a gerbe\footnote{We refer to \cite{2-gerbe} and
  \cite{ima} for the definition of a gerbe, and for additional
  details 
regarding the
  associated 
  cocycle and coboundary equations \eqref{eq:coccases},
  \eqref{eq:cobcases}.}
 on an $S$-scheme $X$. 
For simplicity, in discussing gerbes  we will make two additional
assumptions:

\begin{itemize}
\item  $\pc$ is a $G$-gerbe, for a given $S$-group scheme $G$.

\item $\pc$ is connected.

\end{itemize}

\bigskip

The  first assumption gives us, for any object $x$ in the fibre
category $\pc_U$ above an open set $U \subset X$, an isomorphism of
sheaves on $U$ 
\bee
\label{gger}
\xymatrix{
 G_{|U} \ar[r]^(.37){\sim}& \mathrm{Aut}_{\pc_U}(x)\;.
} 
\end{equation}
 The second assumption asserts that for any pair of objects $x,y \in
 \mathrm{ob}(\pc_U)$ there exists an arrow $x \la y$ in the category 
$\pc_U$. This ensures that the gerbe is described by an element in the
degree 2 \v{C}ech cohomology of $X$ rather than by degree 2 cohomology
with respect to a hypercover of $X$.

\bigskip

Let us   choose  a family
of local objects $x_i \in \pc_{U_i}$, for some open cover $\mathcal{U}
= \coprod_i U_i$ of $X$,  and a family  of arrows
 \begin{equation}
\label{def:lij}
\xymatrix@=15pt{
x_j \ar[rr]^{\phi_{ij}} && x_i 
}\end{equation}
in $\pc_{U_{ij}}$.
 Identifying elements of both  $\mathrm{Aut}_{\pc}(x_i)$   and
 $\mathrm{Aut}_{\pc}(x_j)$ with the corresponding sections of $G$
above $U_i$ and $U_j$, these arrows determine  a family
of  section $\lambda_{ij} \in \Gamma (U_{ij},\,
\mathrm{Aut}(G))$, defined by the commutativity of the diagrams 
\begin{equation}
        \label{def-lamij}
        \xymatrix@R10pt@C=15pt{
   x_j \ar[rr]^{\ga} \ar[dd]_{\phi_{ij}}
    && x_j \ar[dd]^{\phi_{ij}}\\&&\\x_i 
    \ar[rr]_{\lambda_{ij}(\ga)}&& x_i
    }\end{equation} 
\noindent for every  $\ga \in G_{\, |U_{ij}}$. In addition, the arrows $\phi_{ij}$
 determine  a family of
  elements  
 $g_{ijk} \in G_{\mid U_{ijk}}$ for all $(i,j,k)$
by the commutativity of the  diagrams
\begin{equation}
        \label{coc0}
        \xymatrix@R=10pt@C=15pt{
   x_k \ar[rr]^{\phi_{jk}} \ar[dd]_{\phi_{ik}}
    && x_j \ar[dd]^{\phi_{ij}}\\&&\\x_i 
    \ar[rr]_{g_{ijk}}&& x_i
    }\end{equation} 
\noindent above $U_{ijk}$. By conjugation in the sense made clear by
diagram \eqref{def-lamij}, it  follows that the
$\lambda_{ij}$ satisfy the cocycle condition

\begin{equation}
  \label{eq:coc1}
  \lambda_{ij} \, \lambda_{jk} = i(g_{ijk}) \, \lambda_{ik} \,.
\end{equation}
By \cite{ima} lemma 5.1, the $G$-valued cochains $g_{ijk}$
also 
 satisfy the
 cocycle
condition \begin{equation}
  \label{eq:coc2}
  \lambda_{ij}(g_{jkl})\, g_{ijl} = g_{ijk}\, g_{ikl}\,.
\end{equation}
These  two  cocycle equations may be written more compactly as
\bee\label{eq:coccases}
\begin{cases}
\delta^1 \lambda_{ij} & =\,  \: 
i(g_{ijk})   \\
\delta^2_{\lambda_{ij}}(g_{ijk}) & =\, \: 1\;,
\end{cases}
\ee

\medskip

\noi where  $\de^2_{\lambda}$ is 
the $\lambda$-twisted degree 2  \v{C}ech differential
determined by equation \eqref{eq:coc2}. They  may be
jointly viewed as the $(G \la \mathrm{Aut}(G))$-valued 
  \v{C}ech 1-cocycle\footnote{We prefer to emphasize the fact that
   $\lambda_{ij}$ is a 1-cochain since this is more consistent with a
   simplicial definition of the associated cohomology,
 even though it is more customary to
   view   the pair of equations  \eqref{eq:coccases} as a 2-cocycle
   equation, with \eqref{eq:coc1} an auxiliary condition.} equations associated to  the
 gerbe $\pc$,  the open cover $\mathcal{U}$ of $X$, 
and the
 trivializing families of  objects $x_i$ and arrows $\phi_{ij}$ in $\pc$.

\bigskip

Let us  choose   a second
family of local objects $x'_i$  in $\pc_{U_i}$, and of arrows
\begin{equation}
  \label{def:lij1}
\xymatrix@=15pt{
x'_j \ar[rr]^{\phi'_{ij}} && x'_i }
 \end{equation}
above $U_{ij}$. To these correspond a new cocycle pair 
$(\lambda'_{ij},\, g'_{ijk})$. In order to compare this set of arrows  with the
previous one, we
choose (after a harmless  refinement of
the  given
 open cover $\mathcal{U}$ of $X$)  a family of  arrows 
\begin{equation}
\label{def:xi}
\xymatrix@C=30pt{
x_i \ar[r]
^{\chi_i} & x'_i}
\end{equation}
in $\pc_{U_i}$ for all $i$. 
The arrow $\chi_i$ induces by conjugation a section 
   $r_i$ in  the group of sections  $\Gamma (U_i,\,  \mathrm{Aut}(G))$, characterized by the
commutativity of the
 square
\begin{equation}
\label{def-mui}
\xymatrix@R=30pt@C=40pt{ x_i \ar[d]_{\chi_i} \ar[r]^u & x_i 
  \ar[d]^{\chi_i} \\
x'_i  \ar[r]_{r_i(u)} & x'_i
}
\end{equation}
 for all $u \in G$.
The lack of compatibility between these
 arrows $\chi_i$ and the  arrows $\phi_{ij}$, $\phi'_{ij}$
 \eqref{def:lij}, 
 \eqref{def:lij1} is measured by the family of sections $\te_{ij} \in
 \Gamma (U_{ij}
,\, G)$ 
 determined by the commutativity of the
 following diagram:
\begin{equation}
\label{defteij}
\xymatrix@R=15pt@C=40pt{
x_j \ar[r]^{\phi_{ij}} \ar[dd]_{\chi_j} & x_i \ar[d]^{\chi_i} \\
& x'_i \ar[d]^{\te_{ij}} \\
x'_j \ar[r]_{\phi'_{ij}} & x'_i\;.
}
\end{equation} 
 Under the identifications \eqref{gger}, diagram \eqref{defteij}
 induces by conjugation, in a sense made clear by the definition
 \eqref{def-mui} of the auromorphism $r_i$,  
 a commutative diagram of group schemes above ${U_{ij}}$
\[
\label{def:teij1}
\xymatrix@R=15pt@C=40pt{
G \ar[r]^{\lambda_{ij}} \ar[dd]_{r_j} & G \ar[d]^{r_i} \\
& G \ar[d]^{i(\te_{ij})} \\
G \ar[r]_{\lambda'_{ij}} & G\,\,\,,
}
\]
 whose commutativity is expressed by the equation
\begin{equation}
\label{eq:cob1}
\lambda'_{ij} = i(\te_{ij})\, r_i\, \lambda_{ij} \, r_j^{-1}
\end{equation} 
in $\mathrm{Aut}(G)$.

\bigskip

Consider now the diagram\footnote{This diagram  whose 
faces are five
  pentagons and three squares  (as well as those in
  \eqref{diagcoboun-a} and \eqref{diagcoboun-b} below) is  the 1-skeleton of a  Saneblidze-Umble
   cubical model  \cite{S-U}, \cite{loday}
 for the  Stasheff associahedron $K_5$ \cite{sta}.}
\begin{equation}
\label{diagcoboun}
\xymatrix@R=10pt@C=50pt{
&x_k \ar[ddl]_{\phi_{jk}} \ar[rr]^{\phi_{ik}}
\ar@{-}[dd]  
 &&x_i \ar[d]^{\chi_i}  \ar[ddl]_{g_{ijk}}\\
&
&&  x'_i \ar[ddl]^(.7){r_i(g_{ijk}) }
\ar[dddd]^{\te_{ik}} 
\\
   x_j \ar[rr]_(.65){\phi_{ij}}
\ar[dd]_{\chi_{j}}
&\ar@{-}[dd]^{\chi_{k}}
&  x_i
\ar[d]_(.45){\chi_i}
& \\
&
& x'_i \ar[d]^{\te_{ij}} &\\
 x'_j  
\ar[rr]^(.65){\phi'_{ij}}
 \ar[ddd]_{\te_{jk}}
 &\ar[d]
 & x'_i
\ar[ddd]^(.25){ \lambda'_{ij}(\te_{jk})}
 & \\
& x'_k  \ar[ddl]_{\phi'_{jk}} 
\ar@{-}[r]^(.7){\phi'_{ik}}  & \ar[r]& x'_i \ar[ddl]^{g'_{ijk}} \\
&& &\\
  x'_j \ar[rr]_ {\phi'_{ij}} &&  x'_i & \ \ .
}
\end{equation}

\noindent Both the top and the bottom squares commute, since these squares
 are of type  \eqref{coc0}.    So do the back,  the left
 and the top front
vertical squares,
 since all three are of type \eqref{defteij}. The same is true of
 the lower
  front square, and the upper right vertical square, since these two
 are respectively of the
 form \eqref{def-lamij} and \eqref{def-mui}. It follows that the
 remaining lower right  square in the diagram is also
 commutative, since all the arrows in diagram \eqref{diagcoboun} are 
 invertible.  The commutativity
 of this final square is expressed
 algebraically by the equation

\[
\label{eq:cob2}
g'_{ijk} \, \te_{ik}
  = \lambda'_{ij} (\te_{jk})\,\te_{ij} \, r_{i}(g_{ijk})\,.
\]

We say that two   cocycle pairs
   $(\lambda_{ij},\, g_{ijk})$ 
and $(\lambda'_{ij},\,  g'_{ijk})$
 are cohomologous if we are given a pair $(r_i, \te_{ij})$, with
 $r_i \in \Gamma(U_i,\, \mathrm{Aut}(G))$ and $\te_{ij} \in 
\Gamma(U_{ij},\,G)$, satisfying those two equations
\bee
\label{eq:cobcases}
\begin{cases} 
\lambda'_{ij}  &=\: i(\te_{ij})\, r_i\, \lambda_{ij} \, r_j^{-1}\\
g'_{ijk} \, \te_{ik} &
\:= \lambda'_{ij} (\te_{jk})\,\te_{ij} \, r_{i}(g_{ijk})\,.
\end{cases}
\ee
and display this as
\begin{equation}
\label{eq:cob3}
 (\lambda_{ij},\, g_{ijk}) \stackrel{(r_i,\te_{ij})}{\sim}
 (\lambda'_{ij},\,  g'_{ijk})\,.
\end{equation}
The equivalence class of
the cocycle pair $(\lambda_{ij},\,g_{ijk})$ for this relation is
independent of the choices of objects $x_i$ and arrows $\phi_{ij}$ by
from which it was constructed. By definition, it  determines  an element in the  first
non-abelian  \v{C}ech
cohomology set
{\it \v{H}}${}^1(\mathcal{U},\, G \stackrel{i}{\la} \mathrm{Aut}(G))$
with coefficients in the crossed module $i: G \la \mathrm{Aut}(G)$.

\bigskip

\subsection{}
\label{subsection:2-2}

  In \cite{dgg}, the combinatorial description of differential forms
 is used  in order to
define the concepts of connections and curvings on a gerbe.
 For any $S$-group scheme $G$, a
 (relative) connection on a principal $G$-bundle
$P$ above the $S$-scheme
 $X$ may be defined as  a morphism 
\begin{equation}
\label{defcon1} 
\xymatrix@C=30pt{\prb P \ar[r]^\epsilon & \pra P
 }
\end{equation}
between the two pullbacks of $P$ to $\dea$, whose restriction to the
diagonal subscheme 
\[ \Delta: X \hookrightarrow \dea
\]
 is the identity morphism $1_P$.

\bigskip

  This type of  definition of a connection, as a vehicle
  for  parallel
  transport, remains valid for other
structures than principal bundles. In particular,
  for any $X$-group scheme $\Gamma$, a
connection on $\Gamma$ is a morphism of group schemes
 \bee
\label{defmugr}
\mu: \prb \Ga \la
\pra \Ga
\ee
above $\dea$ whose restriction to
  the diagonal subscheme $X \hookrightarrow \dea$
is the identity morphism $1_\Ga$.  When $\Ga$ is the pullback to $X$ of
an $S$-group scheme $G$, the inverse images $\prb G$ and $\pra G$ of
$G_X$ above $\dea$ are canonically isomorphic, so that  the 
connection \eqref{defmugr} is then described by a
$\mathrm{Lie}(\mathrm{Aut}(G))$-valued 1-form  $m$.

\bigskip

 A connection $\mu$ on a group $\Ga$
determines    de Rham differentials 
  \[
\label{def:mucon}
\D^n_{X/S,\,\mu}:\mathrm{Lie}(\Ga) \ot_\osc \Om^n_{X/S} \la
  \mathrm{Lie}(\Ga)  \ot_\osc
  \Om^{n+1}_{X/S}
 \]
(or simply $\D^n_\mu$) defined combinatorially by the formulas
\cite{dgg}
 (A.1.9)-(A.1.11)
and their higher analogues. When $\Ga$ is the pullback of an $S$-group
scheme,
$\D^n_\mu$  is decribed in classical terms as  the deformation 
\eqref{def:dtwist}
\[
\label{def:muconclas}
\D^n_{\mu} :=   \D^n_m 
\]
of the de Rham differential $\D^n$ determined by the associated 1-form $m$.
 When the curvature $\D^1m$ of the connection $\mu$ is trivial, the
 connection is said to be integrable. In that case,  it follows from 
\eqref{dndn1} and \eqref{dndn1a}  that the de Rham differentials
 satisfy  the condition $\D_m^{n+1}\,\D_m^n = 0$ for all $n \neq 1$.

\bigskip 
The curvature 
 of a connection
$\epsilon$  \eqref{defcon1}  on a principal bundle $P$  is the unique arrow
\[\kappa_{\epsilon}:\pra P \la \pra P\]
 such that the following  diagram above 
$\deb$ commutes, with $\epsilon_{ij}$ the pullbacks of $\epsilon$  under the
corresponding projections $p_{ij}: \deb \la \dea$:
\[
        \label{def:kappa1}
        \xymatrix@C=20pt@R=11pt{
   \prc P \ar[rr]^{\epsilon_{12}} \ar[dd]_{\epsilon_{02}}
    && \prb P \ar[dd]^{\epsilon_{01}}\\&&\\\pra 
    P\ar[rr]_{\kappa_{\epsilon}}&& \pra 
    P
    }\] 
  By construction, $\kappa_\epsilon$ is a relative 2-form on $X$ with
  values in the gauge  group $P^{\mathrm{ad}}:=\mathrm{Isom}_G(P,\,P)$ of $P$.

\bigskip

The connection $\epsilon$ on  $P$ induces a connection $\mu_\epsilon$ on 
 the group  $P^{\mathrm{ad}}$, 
  determined by the commutativity of the squares
\begin{equation}
\label{def:ieps}
\:\xymatrix@C=50pt@R=30pt{
p^{\ast}_{1}P \ar[r]^{u} \ar[d]_{\epsilon} & p^{\ast}_{1} P
\ar[d]^{\epsilon}\\ p^{\ast}_{0}P \ar[r]_{\mu_\epsilon (u)} & 
p^{\ast}_{0}P}
\end{equation}
 for all sections $u$ of $\prb(\pa)$. 
By \cite{AK:Bianchi}, \cite{dgg} proposition 1.7, the curvature 2-form $\kappa_\epsilon$
satisfies the Bianchi identity
\begin{equation}
\label{bianchi}
\D^2_{\mu_{\epsilon}}(\kappa_{\epsilon}) = 0.
\end{equation}
For a given family of local sections of $P$, with associated
$G$-valued 
1-cocycles  $g_{ij}$,  the  connection \eqref{defcon1}  is described
by a family of $G$-valued 1-forms $\om_i \in \fg \ot \Om^1_{U_i/S}$, satisfying the gluing
condition 
\bee
\label{con:local}
\om_j = \om_i^{\,\,\ast\,g_{ij}} =  \om_i^{\,g_{ij}} + g_{ij}^{-1}\D g_{ij}
\ee
above $U_{ij}$, for the action of $G$ on $\fg \ot_{\osc}
\Om^1_{U_i/S}$ 
induced by the adjoint right action of
$G$ on $\fg$. A 1-form satisfying
this equation is classically known as a connection form.
The induced curvature $\kappa$ is locally described by the family of 2-forms
\[\kappa_i := \D^1\om_i = \D\om_i + \dfrac{1}{2} [\om_i,\, \om_i], \]
and these satisfy the simpler  \v{C}ech (or  gluing)  condition
\[\kappa_j = \kappa_i^{g_{ij}}\,. \]
Equation \eqref{bianchi} is reflected at the local level  in the
 equation
\[ \D^2_{\,\om_i}\kappa_i = 0\,,\]
which is simply the classical Bianchi identity \eqref{def:bianchiclas}
for the 1-form $\om_i$.

\bigskip

\subsection{}
\label{subsection:2-3}

The notion of a connective structure on a $G$-gerbe $\pc$ is a
categorification of  the notion of  a connection on a principal bundle,
as we will now  recall, following \linebreak \cite{dgg} \S 4.
To $\pc$ is associated its gauge stack $\pac$.  By definition this is
the  monoidal stack $\mathcal{E}q_X(\pc,\,\pc)$ of self-equivalences of the stack $\pc$, the
monoidal structure being defined by the composition of equivalences.
A connection on a $\pc$ is an equivalence between stacks
\begin{equation}
\label{defcon}
\xymatrix{\prb \pc \ar[r]^{\epsilon} & \pra \pc }
 \end{equation}
above $\dea$, together with a natural isomorphism between the  restriction
$\De^{\ast}\epsilon $ of $\epsilon$ to the diagonal subscheme $X$ of $\dea$
 and the identity
morphism $1_{\pc}$. Such a connection $\epsilon$  induces as in
\eqref{def:ieps}
a  connection $\mu$ on the gauge stack $\pac$. 

\newpage

 A curving of $(\pc, \epsilon)$ is a natural
isomorphism $K$ 
 \begin{equation}
         \label{def:kappa2}
         \xymatrix@=13pt{
    \prc \pc \ar[rr]^{\epsilon_{12}} \ar[dd]_{\epsilon_{02}}
     && \prb \pc \ar[dd]^{\epsilon_{01}}
     \ar@{}[dd]_(.55){\, }="1"\\&&\\\pra \pc \ar[rr]_{\kappa}
     \ar@{}[rr]^(.55){\,}="2"&& \pra  \pc \:,
      \ar@{}"1";"2"^(.2){\,}="3"
      \ar@{}"1";"2"^(.8){\,}="4"
      \ar@{=>}"4";"3"^{K}
     }\end{equation}
 for some morphism 
\[
\label{deffake}
\ka:\pra \pc \la \pra \pc
\]
 above $\deb$. It is  determined by the choice of some
explicit  quasi-inverse of the connection $\epsilon$. The arrow $\kappa$ which
 arises
  as part of the definition of $K$ is
called the {\it fake curvature} associated to  the connective structure
$(\epsilon,\, K)$.
It is a global object in the pullback   to 
 $\deb$ of the gauge 
stack $\pac$.

\bigskip

 The connective structure $(\epsilon,\, K)$
determines a 2-arrow  

  \[
\label{defomeg}
\ \ \ \xymatrix@C=35pt{\pra \pc \ar[r]^{\kappa_{013}} \ar[d]_{\kappa_{023}}
&\pra \pc \ar[d]^{\mu_{01}(\kappa_{123})}\ar@{}[d]_(.3){\,}="1"\\
\pra \pc  \ar[r]_{\kappa_{012}}\ar[r]^(.3){\,}="2" &\pra \pc
\ar@{}"1";"2"^(.2){\,}="3"
\ar@{}"1";"2"^(.7){\,}="4"
\ar@{=>}"3";"4"_{\Om}
}
\]

\noi This is  the unique 2-arrow which may be inserted in
diagram
 \begin{equation}
 \label{cube}
         \xymatrix@R=9pt@C=15pt{
    &&\prd \pc\ar[rrrrr]^{\epsilon_{13}} \ar[lldd]_{\epsilon_{03}}
   \ar '[dd][ddddd]^(.3){\epsilon_{23}}
    \ar@{}[rrrrr]_(.2){\,}="1"
    \ar@{}[ddddd]^(.3){\,}="2"
     &&&&& \prb \pc
    \ar[lldd]_(.6){\epsilon_{01}}
    \ar@{}[lldd]^(.7){\,}="6"
 \ar@{}[lldd]^(.6){\,}="21"
 \ar[ddddd]^{\kappa_{123}}
 \ar@{}[ddddd]^(.6){\,}="22"\\
    &&&&&&&\\
    \pra \pc\ar[rrrrr]^(.6){\kappa_{013}}
     \ar@{}[rrrrr]^(.7){\,}="5"
    \ar[ddddd]_{\kappa_{023}}
    \ar@{}[ddddd]^(.5){\,}="9"
    &&&&&\pra \pc \ar[ddddd]^(.4){\mu_{01}(\kappa_{123})} 
    \ar@{}[ddddd]_(.8){\,}="17"
    &&\\
    &&&&&&&\\
    &&&&&&&\\
    &&\prc \pc \ar[lldd]^(.4){\epsilon_{02}} 
    \ar@{}[lldd]^(.4){\,}="13"
    \ar@{}[lldd]^(.6){\,}="10"
    \ar '[rrr]^(.7){\epsilon_{12}} [rrrrr]
    \ar@{}[rrrrr]_(.3){\,}="14"
    &&
    &&&
    \prb \pc \ar[lldd]^{\epsilon_{01}}
    \\
    &&&&&&&
    \\ \pra \pc \ar[rrrrr]_{\kappa_{012}}
     \ar@{}[rrrrr]^(.8){\,}="18"
     &&&&& \pra \pc && \qquad .
    \ar@{}"1";"2"^(.3){\,}="3"
     \ar@{}"1";"2"^(.8){\,}="4"
   \ar@{=>}"3";"4"^{K_{123}}
   \ar@{}"5";"6"^(.4){\,}="7"
   \ar@{}"5";"6"^(.75){\,}="8"
   \ar@{}"9";"10"^(.3){\,}="11"
  \ar@{}"9";"10"^(.7){\,}="12"
   \ar@{=>}"7";"8"^{K_{013}}
   \ar@{=>}"11";"12"^{K_{023}}
    \ar@{}"13";"14"^(.3){\,}="15"
     \ar@{}"13";"14"^(.6){\,}="16"
    \ar@{=>}"15";"16"_{K_{012}}
     \ar@{}"17";"18"^(.2){\,}="19"
       \ar@{}"17";"18"^(.8){\,}="20"
        \ar@{=>}"19";"20"_{\Omega}
        \ar@{}"21";"22"^(.35){\,}="23"
\ar@{}"21";"22"^(.25){\,}="25"
        \ar@{}"21";"22"^(.8){\,}="24"
 \ar@{}"21";"22"^(.7){\,}="26"
  \ar@{}"23";"24"_{M_{01}(\kappa_{123})}
\ar@{=>}"25";"26"
}
\end{equation}
 so that the
two composite  2-arrows  
\[\xymatrix@C=90pt{
\prd \pc \ar@/^1pc/[r]^{
\mu_{01}(\kappa_{123})  \, \kappa_{013}\,\eps_{03}} \ar@/^1pc/[r]_{\,}="1"
\ar@/_1pc/[r]_{\eps_{01}\, \eps_{12} \, \eps_{23} }
\ar@/_1pc/[r]^{\,}="2"
& \pra \pc
\ar@{}"1";"2"^(.1){\,}="3"
\ar@{}"1";"2"^(.9){\,}="4"
\ar@{=>}"3";"4"
}\] 
which  may  be constructed by composition of 2-arrows in \eqref{cube}
coincide.

\bigskip

\noi  This 2-arrow $\Om$ may also   be viewed as a  1-arrow above $\dec$ in the
 gauge group $\pac$,  or even as an arrow in the stack
 $\mathrm{Lie}(\pac) \ot_\osc \Om^3_{X/S}$ of relative 
 $\mathrm{Lie}(\pac)$-valued 3-forms on $X$.
 Returning
 to the combinatorial definition \cite{dgg} (A.1.10) of the de Rham 
differential, we may finally  view $\Om$, by 
horizontal composition with appropriate 1-arrows, as a 
  1-arrow in $\pac$  whose source
object is the identity arrow $I_{\pc^{\mathrm{ad}}}$:
\begin{equation}
\label{tom}
\xymatrix@C=30pt{
I \ar[r]^(.38){\Om} & \D^2_{\mu}(\kappa^{-1}).
}
\end{equation}
Denoting  the twisted differential $\D^2_{\mu}$ by the expression 
$\D  + [\mu,\,\ ]$ to which it reduces when appropriate
trivializations have been chosen,   the 3-curvature arrow $\Om$ \eqref{tom}
is  described by the equation \eqref{eqom1}. By 
\cite{dgg} theorem 4.4 it satisfies another
relation, described by the  cubical pasting diagram  \cite{dgg}
(4.1.24), and which may be expressed by the higher Bianchi
 identity\footnote{See  \cite{dgg} (4.1.28) for a proof of this
identity.}
\eqref{eqom2}. The pair of equations \eqref{eqom1} and \eqref{eqom2}
may now be thought of as a categorified version, satisfied by the pair
of $\pac$-valued forms $(\kappa, \, \Om)$, of the classical
Bianchi identity \eqref{bianchi}, and can  be written in symbolic form as
\[ \D^2_{\mu,\,\mathcal{K}} (\kappa,\,\Om) = 0\,,\]
 where   $ \D^n_{\mu,\,\mathcal{K}} $     is 
 the twisted de Rham
  differential on $\mathrm{Lie} (\pac)$-valued $n$-forms determined by
   twisting data $(\mu,\, \kc)$ associated to the given connective
   structure on $\pc$.

\section{\v{C}ech-de Rham cocycles}
\label{section3}

\subsection{}
\label{subsection:3-1}

We observed    in  section \ref{section2}.1 that  a gerbe could be
expressed in cocyclic terms, once local trivializations were
chosen. We 
will now show that this is also the case for the  connection
$\epsilon$.
 We choose, for each $i \in I$, an arrow
 \begin{equation}
\label{def:gai}
\gamma_i:\eps  \prb x_i \la \pra x_i
\end{equation}
 in $\pra \pc_{U_i}$ such that $\Delta^{\ast}\ga_i = 1_{x_i}$.
The arrow $\gamma_i$ determines by conjugation a connection
\[m_i:\prb G_{|U_i} \la \pra G_{|U_i} \]
 on
the pullback  $ G_{|U_i}$ of the  group
$G$ above the open set $U_i \subset X$. The arrow $m_i$ is  described,
for any section  $g \in \Ga(\dea_{U_i},\, \prb G)$,
by
the commutativity of the  diagram
\bee
\label{defconmi}
\xymatrix@C=40pt@R=20pt{
\epsilon \prb x_i \ar[r]^{\epsilon  (g)} \ar[d]_{\ga_i}& 
\epsilon \prb x_i \ar[d]^{\ga_i} \\\pra x_i \ar[r]_{m_i(g)} & \pra x_i\;.
}
\ee

\bigskip

The pair $(\phi_{ij},\, \ga_i)$ determines 
 a family of arrows $\gamma_{ij}$ in  the pullback
 $G_{\Delta^1_{U_{ij}}}$ of $G$, defined by the
 commutativity of the diagram 
\begin{equation}
\label{diag:gaij}
\xymatrix@C=50pt@R=15pt{
\eps\prb x_j \ar[r]^{\gamma_j}
\ar[dd]_{\eps \prb \phi_{ij}} 
 & \pra x_j \ar[d]^{\pra \phi_{ij}}\\
& \pra x_i \ar[d]^{\gamma_{ij}}\\
\eps \prb x_i \ar[r]_{\gamma_i} & \pra x_i 
}
\end{equation}
By conjugation, this  determines
a commutative diagram
\begin{equation}
\label{diag:gaij1}
\xymatrix@C=50pt@R=15pt{
\prb G \ar[r]^{m_j}
\ar[dd]_{ \prb \lambda_{ij}} 
 & \pra  G \ar[d]^{\pra \lambda_{ij} } \\
& \pra  G       \ar[d]^{i(\ga_{ij})} \\
 \prb G \ar[r]_{m_i} & \pra G 
}
\end{equation}
so  that  the equation 
\begin{equation}
    \label{cocep13}
     i(\ga_{ij})\, (\pra \lambda_{ij})\,m_{j} \, (\prb 
     \lambda_{ij})^{-1}=  m_{i}\:.
        \end{equation}
of  \cite{dgg} (6.1.2)  is satisfied.

\bigskip

We may restate  (\ref{cocep13}) as
\bee
\label{cocep13new}
 i(\ga_{ij})\,\, [(\pra \lambda_{ij})\,m_{j}\, (\pra \lambda_{ij})^{-1}]
= m_i \,\,[\prb 
     \lambda_{ij} \,\,(\pra \lambda_{ij}^{-1})]\,,
\ee
an equation all  of whose factors are $\mathrm{Aut}(G)$-valued 1-forms
on $U_{ij}$ and therefore commute with each other. 
In  the notation introduced in   \eqref{cocep13clas-b},
equation \eqref{cocep13new} can be rewritten as   
\bee
\label{cocep13clas0}
 {}^{\lambda_{ij}\,\ast\,}\!m_j = m_i - i(\ga_{ij})\,,
\ee
 or more classically as
\bee
\label{cocep13clas}
 {}^{\lambda_{ij}\,}\!m_j = m_i -
\lambda_{ij}\,\,\D\lambda_{ij}^{-1} -i(\ga_{ij})\,.
\ee
This is is the analogue for the  $\mathrm{Aut}(G)$-valued forms $m_i$
and $\lam_{ij}$
of the classical expression (\ref{con:local}) for a connection form,
but now categorified by  the insertion of an  additional summand 
 $- i(\ga_{ij})$.

\newpage

Consider now the following diagr
in $\pc_{\Delta^1_{U_{ijk}}}$:
\begin{equation}
\label{diagcoboun-a}
\xymatrix@R=10pt@C=50pt{
 & \epsilon \prb x_k \ar[ddl]_{\epsilon \prb \phi_{jk}}
\ar[rr]^{\epsilon \prb \phi_{ik}}
\ar@{-}[dd]
  &&\epsilon \prb x_i \ar[d]^{\gamma_i} 
 \ar[ddl]_{\epsilon \prb g_{ijk}}\\
&
&&  \pra x_i \ar[ddl]^(.7){m_i(\prb g_{ijk}) } \ar@{<-}[dddd]^{\gamma_{ik}} \\
 \epsilon \prb  x_j \ar[rr]_(.65){\epsilon \prb \phi_{ij}}
\ar[dd]_{\gamma_{j}}
& \ar@{-}[dd]& \epsilon \prb  x_i
\ar[d]_(.45){\gamma_i}
& \\
&& \pra x_i \ar@{<-}[d]^{\gamma_{ij}} &\\
\pra  x_j  
\ar[rr]^(.65){\pra \phi_{ij}}
 \ar@{<-}[ddd]_{\ga_{jk}}
 &\ar[d]_{\gamma_{k}}  & \pra x_i
\ar@{<-}[ddd]^(.25){ \lambda_{ij}(\ga_{jk})}
 & \\
& \pra x_k  \ar[ddl]_{\pra \phi_{jk}} 
\ar@{-}[r]^(.7){\pra \phi_{ik}}
 &\ar[r]& \pra x_i \ar[ddl]^{\pra g_{ijk}} \\
&& &\\
  \pra x_j \ar[rr]_ {\pra \phi_{ij}} && \pra  x_i   & 
}
\end{equation}
Of the eight faces of this cube, seven  are known to be commutative.
 It follows that the remaining lower square on the right vertical side
 is also commutative. This is the  square
\bee
\label{eq:gaij11}
\xymatrix@C=50pt@R=20pt{
\pra x_i \ar[r]^{\pra g_{ijk}} \ar[dd]^{\ga_{ik}} & \pra x_i 
\ar[d]^{\lambda_{ij}(\ga_{jk})}\\
& \pra x_i \ar[d]^{\ga_{ij}}\\
\pra {x_i} \ar[r]_{m_i(\prb g_{ijk})} & \pra x_i\,, 
}
\ee
whose  commutativity   corresponds to the equation
\[
\label{cocep5}
\ga_{ij}\, \,(\pra \lambda_{ij} (\gamma_{jk}))
\:= \:  m_{i}(\prb g_{ijk})\,\, \ga_{ik}\,\,\,(\pra g_{ijk})^{-1}
\]
in other words to
 the equation \cite{dgg} (6.1.7), all of  whose factors are
 $G$-valued 1-forms on $U_{ijk}$.
We may rewrite this as 
\[ \ga_{ij}\,\, \pra \lam_{ij}(\ga_{jk}) = (m_i(\prb g_{ijk})\,\pra
g_{ijk}^{-1})\,
(\pra g_{ijk}\, \ga_{ik}\, \pra g^{-1}_{ijk})
\]
so that, taking into account the equation \eqref{eq:coc1}, we finally
obtain (in additive notation)
\[
\label{cocep5clas}
\ga_{ij} + \lam_{ij}(\ga_{jk}) - \lam_{ij}\lam_{jk}(\lam^{-1}_{ik}(\ga_{ik}))
  =  dg_{ijk}\,g^{-1}_{ijk} + [m_i,\, g_{ijk}]\,,
\]
with bracket defined by \eqref{def:bra1}
an equation which can be written in abbreviated form as
\bee
\label{cocep5clas1}
\de^1_{\lambda_{ij}}(\ga_{ij}) = \D_{m_i}g_{ijk}\,\,g^{-1}_{ijk} \,.
\ee

\bigskip

\subsection{}
\label{subsection:3-2}

We  now describe in similar terms  the  curving $K$ and   the  fake curvature
$\kappa$  of diagram \eqref{def:kappa2}. 
  Just as  we associated to the connection  $\eps$ \eqref{defcon} a
  family
of arrows $\ga_i$ \eqref{def:gai}, we now  choose, for each $i \in I$, an arrow
\begin{equation}
\label{def:dei}
\xymatrix{ \kappa \pra x_i \ar[r]^{\de_i} &\pra x_i }
\end{equation}
 in the category  $\pc_{\Delta^{2}_{U_i}}$, whose restriction to the degenerate
 subsimplex $s\Delta^{2}_{U_i}$ of $\Delta^{2}_{U_i}$ is the identity.
To the curving $K$  is associated  a family of
``$B$-field''  $\fg$-valued 2-forms $B_i \in \fg \ot \Om^2_{U_i}$, 
characterized by the commutativity of the following
diagram\footnote{The chosen  orientation of the arrow $B_i$ 
is consistent with that in \cite{dgg}.} in which
an expression such as $\ga_i^{12}$ is the pullback of $ \ga_i$ by the
corresponding projection $p_{12}: \deb \la \dea$:

\begin{equation}
\label{def:bi}
\xymatrix@C=40pt@R=25pt{
\eps_{01}\eps_{12}(\prc x_i)
\ar[d]_{\eps_{01}\ga_i^{12}}
 \ar[r]^{K(\prc x_i)} & \kappa \eps_{02} (\prc x_i) \ar[d]^{\kappa \ga_i^{02}}\\
\eps_{01}(\prb x_i) \ar[d]_{\ga_i^{01}} &\kappa  \pra x_i \ar[d]^{\de_i}\\
\pra x_i \ar@{<-}[r]_{B_i}  & \pra x_i
}
\end{equation}

\noi Let us now  define a family of $G$-valued 2-forms $\nu_i$ on $U_i$  by  the
equations
\begin{equation}
   \label{ifi}
  \nu_{i} :=  \D^1m_{i} -  i(B_{i}) 
    \end{equation}
in $ \mathrm{Lie
\:  Aut}(G) \ot \Om^2_{U_i}$, in other words by the
 commutativity of the diagram
\bee
\label{biconj}
\xymatrix@C=40pt@R=20pt{
\prc G \ar@{=}[r] \ar[d]_{m_i^{12}} & \prc G \ar[d]^{m_i^{02}} \\
\prb G \ar[d]_{m_i^{01}} & \pra G \ar[d]^{\nu_i}\\
\pra G_i  \ar@{<-}[r]_{i(B_i)} & \pra G \;.
}
\ee  
 By comparing diagram \eqref{biconj} with the
    conjugate of diagram \eqref{def:bi}, we  see that $\nu_i$ is
   simply the conjugate of the arrow $\de_i$.
It can    therefore  described  by
    the commutativity of the diagram 
\begin{equation}
\label{defnui}
\xymatrix@C=40pt@R=20pt{
\kappa \pra x_i \ar[r]^{\kappa (g)} \ar[d]_{\de_i}& 
\kappa \pra x_i \ar[d]^{\de_i} \\\pra x_i \ar[r]_{\nu_i(g)} & \pra x_i
}
\end{equation}
for all $g \in \Ga(\Delta^2_{U_i/S},\,\pra  G)$, just as the
connection  $m_i$ was described  by
diagram \eqref{defconmi}.

\bigskip

 We also define a family of 2-forms $\de_{ij}$
 by the commutativity of
 the diagram
\bee
\label{def:deij1}
\xymatrix@R=20pt@C=40pt{
\pra x_i \ar[r]^{\lambda_{ij}(B_j)}  \ar[d]_{\de_{ij}}& \pra x_i 
\ar[d]^{\ga_{ij}^{01}}\\
\pra x_i  \ar[d]_{\ga_{ij}^{02}}
 & \pra x_i  \ar[d]^{m_i^{01}(\ga_{ij}^{12})}\\
\pra x_i \ar[r]_{B_i}& \pra x_i\,,
}
\ee
i.e., since all terms commute, by the equation
\[
\label{eq:bij1}
 \de_{ij} := \lambda_{ij}(B_j) -  B_i - \D^1_{m_i}(- \ga_{ij})
\]
in $\mathrm{Lie}(G) \ot \Om^2_{U_i/S}$.
In \v{C}ech-de Rham notation, this is 
\begin{equation}
\label{eq:bij1a}
 \de_{ij} :=\de^0_{\lambda_{ij}}(B_i) -\D^1_{m_i}(- \ga_{ij})\;,
\end{equation}
and in   classical notation
\[
\label{eq:bij11}
 \de_{ij}:=  \lambda_{ij}(B_j) -  B_i + \D\ga_{ij} -
\dfrac{1}{2}[\ga_{ij},\, \ga_{ij}] + [m_i,\, \ga_{ij}]\;.
    \]

\bigskip

\noi Here is another  characterization of $\de_{ij}$:

\bigskip

\begin{lemma}
  For every pair $(i,j) \in I$, the
  analogue 
\begin{equation}
\label{diag:deij}
\xymatrix@C=50pt@R=15pt{
\kappa\pra x_j \ar[r]^{\delta_j}
\ar[dd]_{\kappa \pra \phi_{ij}} 
 & \pra x_j \ar[d]^{\pra \phi_{ij}}\\
& \pra x_i \ar[d]^{\delta_{ij}}\\
\kappa \pra x_i \ar[r]_{\de_i} & \pra x_i \,.
}
\end{equation}
of diagram  \eqref{diag:gaij} is  commutative.
\end{lemma}

\newpage

\noi {\bf Proof}: Consider the diagram

\begin{equation}
\label{pro2}
\xymatrix@R=25pt@C=18pt{
& \kappa \eps_{02}(\kappa\prc x_j)
\ar@{-}[d] 
\ar[rrr]^{\kappa \ga_j^{02}}
&&&\kappa  \pra x_j  
 \ar[dd]^(.45){\kappa\pra \phi_{ij}} 
 \ar[r] ^{\de_j} & \pra x_j \ar[d]_{\pra \phi_{ij}}
   \\
\eps_{01} \eps_{12} (\prc x_j)
\ar[ddddd]^{\eps_{01}\eps_{12}(\prc \phi_{ij})}
 \ar[rr]^(.7){\ga_j^{12}}\ar[ur]^{K(\prc x_j)} \
&\ar[dddd]^(.6){\kappa \eps_{02}(\prc \phi_{ij})}
& \eps_{01}( \prb x_j) \ar[r]^(.6){\ga_j^{01}}
 \ar[ddd]_(.4){\eps_{01}(\prb \phi_{ij})}
 & \pra x_j \ar@{<-}[urr]^(.3){B_j}
\ar[d]_{\pra \phi_{ij}}   &
&   \pra x_i  \ar[d]_{\de_{ij}} \\
& &&   \pra x_i   \ar@{<-}[urr]^(.35){\lambda_{ij}(B_j)} 
  & \kappa \pra x_i\ar[r]^{\de_i} & \pra x_i  
 \ar[ddd]_{\nu_i(\ga_{ij}^{02})} 
\\
&&& &&\\
&&\eps_{01}(\prb x_i) \ar[r]_{\ga_i^{01}}  & \pra x_i
\ar@{<-}[uu]_{\ga_{ij}^{01}}&\\
&\kappa  \eps_{02}(\prc x_i) 
\ar@{-}[r]
 &\ar@{-}[r]^{\kappa\ga_i^{02}} &\ar[r]&
\kappa \pra x_i \ar[r]^{\de_i}
\ar@{<-}[uuu]^{\kappa \ga_{ij}^{02}}& \pra x_i\\
\eps_{01}\eps_{12}(\prc x_i) \ar[ur]^{K(\prc x_i)}
\ar[rr]_{\eps_{01}(\ga_i^{12})}  && \eps_{01}(\prb x_i)
\ar@{<-}[uu]^(.3){\eps_{01}(\ga_{ij}^{12})}\ar[r]_{\ga_i^{01}}& \pra x_i
\ar@{<-}[uu]^(.3){m_{i}^{01}(\ga_{ij}^{12})} \ar@{<-}[urr]_{B_i}&&\:.
}
\end{equation}

\noi  Diagrams  \eqref{def:bi}, \eqref{def:deij1} and \eqref{defnui}
imply that  all
squares in \eqref{pro2}  are commutative\footnote{This is true for
diagram \eqref{def:deij1}  since $\nu_i(\ga_{ij}^{02}) =
\ga_{ij}^{02}$. }, except possibly the rear right
upper one. This remaining square (\ref{diag:deij}) 
is therefore also commutative. \hspace{6.9cm}$\Box$

\bigskip

\noi Conjugating diagram
\eqref{diag:deij}, we obtain as in \eqref{diag:gaij1}  a square 
\[
\xymatrix@C=50pt@R=15pt{
 \pra G \ar[r]^{\nu_j}
\ar[dd]_{ \pra \lambda_{ij}} 
 & \pra  G \ar[d]^{\pra \lambda_{ij} } \\
& \pra  G       \ar[d]^{i_{\de_{ij}}} \\
\kappa \pra G \ar[r]_{\nu_i} & \pra G\,, 
}
\]
 whose  commutativity  is expressed algebraically as
\bee \label{cockap01}
i(\de_{ij})\,  (\pra \lambda_{ij})\,\nu_{j} \:=\:  \nu_{i}\,
        (\pra \lambda_{ij})\,.
        \end{equation}
\newpage

\noi In additive notation, this is  equation 
\begin{equation}
    \label{cockap1}
    {}^{\lambda_{ij}\,}\!\nu_{j} = \nu_{i} - i(\delta_{ij})\:,
    \end{equation}
in other words
\[
\label{cockap1-cech}
\de^{0}_{\lambda_{ij}}\nu_i = -\, i(\delta_{ij})\:.
    \]
It is instructive to note that  this equation can be derived directly
 from
equation \eqref{cocep13clas} and the definitions 
\eqref{ifi} and \eqref{eq:bij1a}  of $\nu_i$ and $\de_{ij}$.
First of all, observe  that by  \eqref{funct:d1}
\bee
\label{eq:d1mj}
\D^1({}^{\lambda_{ij}\ast\,}\!m_i) = {}^{\lambda_{ij}\,}\!(\D^1m_i) \;.
\ee
One then computes
\begin{eqnarray*}
{}^{\lambda_{ij}\,}\!\nu_j &=&{}^{\lambda_{ij}\,}\!(\D^1(m_j) - i_{B_j})\\
&=&\D^1({}^{\lambda_{ij}\ast\,}\!m_j) - i(\lambda_{ij}(B_j))\\
&=& \D^1(m_i - i(\ga_{ij})) - i(B_i + \D^1_{m_i}(-\ga_{ij}) +
  \de_{ij})\\
&=& \D^1m_i - \D^1( i(\ga_{ij}))  - [m_i, \ga_{ij}] - i(B_i)
- i( \D^1m_i(-\ga_{ij})) -
 i(\de_{ij})\;.
\end{eqnarray*}
Since the homomorphism $i$ commutes with $\D^1m$ and 
$ [m_i,\,i(\ga_{ij})] = i([m_i,\,\gamma_{ij}])$,  the summands
$i(\D^1m(-\ga_{ij}))$ and $\D^1(i(\ga_{ij})) + [m_i,\,\ga_{ij}]$ cancel
out. The first two remaining summands describe $\nu_i$,
so 
that equation  \eqref{cockap1} is satisfied.

\bigskip

 In the same vein, the analogue   for the fake curvature $\kappa$ of 
\eqref{eq:gaij11} is the following assertion.

\bigskip

\begin{lemma}
The diagram
\bee
\label{eq:dij3}
\xymatrix@C=50pt@R=20pt{
\pra x_i \ar[r]^{\pra g_{ijk}} \ar[dd]^{\delta_{ik}} & \pra x_i 
\ar[d]^{\lambda_{ij}(\delta_{jk})}\\
& \pra x_i \ar[d]^{\delta_{ij}}\\
\pra {x_i} \ar[r]_{\nu_i(\pra g_{ijk})} & \pra x_i 
}
\ee
 is commutative.
\end{lemma}

\newpage

\noi {\bf Proof:} By \eqref{diag:deij}, \eqref{coc0} and
\eqref{defnui},
all squares in the  diagram 
\begin{equation}
\label{diagcoboun-b}
\xymatrix@R=15pt@C=50pt{
 & \kappa \pra x_k \ar[ddl]_{\kappa \pra \phi_{jk}}
\ar[rr]^{\kappa \pra \phi_{ik}}
\ar@{-}[dd]
  &&\kappa \pra x_i \ar[d]^{\delta_i} 
 \ar[ddl]_{\kappa \pra g_{ijk}}\\
&
&&  \pra x_i \ar[ddl]^(.7){\nu_i(\pra g_{ijk}) } \ar@{<-}[dddd]^{\delta_{ik}} \\
 \kappa \pra  x_j \ar[rr]_(.65){\kappa \pra \phi_{ij}}
\ar[dd]_{\delta_{j}}
& \ar@{-}[dd]_{\delta_{k}}& \kappa \pra  x_i
\ar[d]_(.45){\delta_i}
& \\
&& \pra x_i \ar@{<-}[d]^{\delta_{ij}} &\\
\pra  x_j  
\ar[rr]^(.65){\pra \phi_{ij}}
 \ar@{<-}[ddd]_{\delta_{jk}}
 &\ar[d]  & \pra x_i
\ar@{<-}[ddd]^(.25){ \lambda_{ij}(\delta_{jk})}
 & \\
& \pra x_k  \ar[ddl]_{\pra \phi_{jk}} 
\ar@{-}[r]^(.6){\pra \phi_{ik}} & \ar[r]& \pra x_i \ar[ddl]^{\pra g_{ijk}} \\
&& &\\
  \pra x_j \ar[rr]_ {\pra \phi_{ij}} && \pra  x_i & 
}
\end{equation}
are commutative, except possibly the lower right-hand one. It follows
that the latter one,  which is
simply \eqref{eq:dij3}, also commutes.
\qed

\bigskip
 
 The commutativity  of \eqref{eq:dij3}  corresponds to  equation
\[
\de_{ij}\, \,(\pra \lambda_{ij})(\de_{jk})
\:= \:  \nu_{i}(\pra g_{ijk})\,\, \de_{ik}\, (\pra g_{ijk})^{-1}\:,
\]
an equation whose terms are
 $G$-valued 2-forms on $U_{ijk}$. By the same reasoning as for
  \eqref{cocep5clas1}, this can be written additively as
\[\de_{ij} + \lam_{ij}(\de_{jk}) 
- \lam_{ij}\lam_{jk}(\lam^{-1}_{ik}(\de_{ik}))
= [\nu_i, \, g_{ijk}]\,,\]
or, in the compact form of \cite{dgg} (6.1.15), as
\begin{equation}
    \label{cockap2}
    \de^{1}_{\lambda_{ij}}(\delta_{ij}) = [\nu_{i},\,g_{ijk}]\, .
    \end{equation}  
Just we were able to derive (\ref{cockap1}) directly  from (\ref{cocep13clas})
and the definitions (\ref{ifi})  and (\ref{eq:bij1a}), we now show
that  it is possible to deduce
(\ref{cockap2}) from (\ref{eq:bij1a}),(\ref{ifi})  and
(\ref{cocep5clas1}). First of all,

\begin{eqnarray}
\de^1_{\lambda_{ij}}(\de_{ij}) &=& \notag
 \de^1_{\lambda_{ij}}( \de^0_{\lambda_{ij}}(B_i) - \D^1_{m_i}(-
 \ga_{ij})  ) \\
&=& \de^1_{\lambda_{ij}}\de^0_{\lambda_{ij}}(B_i) 
- \de^1_{\lam_{ij}}\D^1_{m_i}(-\gamma_{ij})\;. \label{cockap2a} 
\end{eqnarray}

We now wish to assert that the \v{C}ech differential  $\de^1_{\lam_{ij}}$
and de Rham differential $\D^1_{m_i}$ in \eqref{cockap2a} commute with
each other, despite the fact that the 1-form $\gamma_{ij}$ takes its
values in a non-commutative group $G$, and that $\D^1_{m_i}$ is not a
homomorphism. For this we simplify our notation, by setting
\bee
\label{gatilde}
\widetilde{\ga}_{ij}:= - \gamma_{ij} \in \fg \ot \Om^1_{U_{ij}}\ee
 and 
\begin{equation*}
\lam_{ijk}:= \lam_{ij}\,\lam_{jk}\, \lam^{-1}_{ik} \in
\Gamma(U_{ijk},\, \mathrm{Aut}(G_i)).
\end{equation*}

\noi Equation  \eqref{cocep5clas1} can be restated as

\begin{equation}
\label{cocep5clas-b}
\de^1_{\lam_{ij}}{\wtg}:=
 \widetilde{\gamma}_{ij}+ \lam_{ij}(\widetilde{\ga}_{jk}) - \lam_{ijk}
 (\wtg_{ik})= -\, \D g_{ijk}\, g_{ijk}^{-1} - [m_i,\, g_{ijk}]\;.
\end{equation}

\begin{lemma} The following equality between $G$-valued 2-forms above
  $U_{ijk}$ is satisfied:
\bee
\label{comd1} 
\D^1_{m_i}\de^1_{\lam_{ij}}(\wtg_{ij}) =
\de^1_{\lam_{ij}}\D^1_{m_i}(\wtg_{ij})\;.
\end{equation}
\end{lemma}

\noi {\bf Proof:} We compute the left-hand  side of the equation
\eqref{comd1}, taking into account
  the quadraticity  equation
\eqref{d1add}
\begin{eqnarray*}
\label{eq:d1del}
\D^1_{m_i}\,\de^1_{\lam_{ij}}(\wtg_{ij}) & =& 
\, \D_{m_i}(\wtg_{ij})
+ \D^1_{m_i}(\lam_{ij}(\wtg_{jk})) +
\D^1_{m_i}(-\lam_{ijk}(\wtg_{ik})) +  \notag \\
         &  & \hspace{4cm}  
+\,  [\wtg_{ij},\, \lam_{ij}(\wtg_{jk})] -  [\wtg_{ij},\,
\lam_{ijk}(\wtg_{ik})]
- \notag \\
&& \hspace{6cm} - \, [\lam_{ij}(\wtg_{jk}),\, \lam_{ijk}(\wtg_{ik})]
\notag  
 \\
&=&  \,  \D_{m_i}(\wtg_{ij})
+ \D^1_{m_i}(\lam_{ij}(\wtg_{jk}))   
   -   \D^1_{m_i}(\lam_{ijk}(\wtg_{ik})) + \notag  \\
&&  \hspace{3cm} 
 +\, [ \lam_{ijk}(\wtg_{ik}),\,   \lam_{ijk}(\wtg_{ik})]      
 +\, [\wtg_{ij},\, \lam_{ij}(\wtg_{jk})] -  \notag  \\
&& \hspace{5cm}  -\,  [\wtg_{ij}
 + \lam_{ij}(\wtg_{jk}),\, \lam_{ijk}(\wtg_{ik})] \;. \notag 
\end{eqnarray*}
\noi We now compute the right-hand side of \eqref{comd1}:
\begin{eqnarray}
\label{eq:del1d}
\de^1_{\lam_{ij}}\,\D^1_{m_i}(\wtg_{ij}) &=& \D^1_{m_i}(\wtg_{ij}) +
\lam_{ij}(\D^1_{m_j}(\wtg_{jk})) - \lam_{ijk}(\D^1_{m_i}(\wtg_{ik}))\;.
\end{eqnarray}

\noi By  \eqref{cocep13clas0} and by the functoriality property 
\eqref{funct:d1}, we find that 
\begin{align*}
\lam_{ij}
(\D^1_{m_j}(\wtg_{jk})) &= \D^1_{{}^{\,\,\lam_{ij}\,\ast\,}\!m_j}
(\lam_{ij}(\wtg_{jk}))\\
&= \D^1_{m_i}(\lam_{ij}(\wtg_{jk})) + 
\, [\wtg_{ij},\, \lam_{ij}(\wtg_{jk})]\\ \intertext{and by \eqref{funct:dnm1}}
\lam_{ijk}(\D^1_{m_i}(\wtg_{ik})) &= 
\D^1_{\,\,{}^{\lam_{ijk}\ast\,}\!m_i}(\lam_{ijk}(\wtg_{ik})) \\
&= \D^1_{m_i}(\lam_{ijk}(\wtg_{ik}))  
 +\, [[\lam_{ijk},\,m_i],\,\lam_{ijk}(\wtg_{ik})]  +\,\\&
\hspace{4.5cm} +\,
[ \lam_{ijk} \,\D\lam^{-1}_{ijk},\,\,\lam_{ijk}(\wtg_{ik})]\;.
\end{align*}
Inserting these expressions for $\lam_{ij}
(\D^1_{m_j}(\wtg_{jk}))$ and $\lam_{ijk}(\D^1_{m_i}(\wtg_{ik}))$ into
the right-hand side of \eqref{eq:del1d}
we find the following expression for
$\de^1_{\lam_{ij}}\,\D^1_{m_i}(\wtg_{ij})$:
\begin{align*}
\de^1_{\lam_{ij}}\,\D^1_{m_i}&(\wtg_{ij}) =
\D^1_{m_i}(\wtg_{ij}) +  \D^1_{m_i}(\lam_{ij}\wtg_{jk})\,
 +\,[\wtg_{ij},\, \lam_{ij}(\wtg_{jk})] \,- \\ & \hspace{2mm}
 - 
 \,\D^1_{m_i}(\lam_{ijk}(\wtg_{ik})
\, - \,  [[\lam_{ijk},\,m_i],\,\lam_{ijk}(\wtg_{ik})] \,-\,
[ \lam_{ijk} \,\D\lam^{-1}_{ijk},\,\lam_{ijk}(\wtg_{ik})] \,- \\&
\hspace{4cm}  
\,  - \,  \D^1_{m_i}(\lam_{ijk})(\wtg_{ik}) - \,\,
[ \lam_{ijk} \,\D\lam^{-1}_{ijk},\,\lam_{ijk}(\wtg_{ik})]\,.
\end{align*}
\noi Comparing this  with the expression  \eqref{eq:d1del} for
$  \D^1_{m_i}\,\de^1_{\lam_{ij}}(\wtg_{ij})$, we see that the
 equation \eqref{comd1} is satisfied if and only if
\begin{multline*}
[\wtg_{ij} + \lam_{ij}(\wtg_{jk}) - \lam_{ijk}(\wtg_{ik}),\,
\lam_{ijk}(\wtg_{ik})]  =
 [[\lam_{ijk},\, m_i],\, \lam_{ijk}(\wtg_{ik})] + \\ +\,
 [\lam_{ijk}\,\D\lam^{-1}_{ijk},\, \lam_{ijk}(\wtg_{ik})]\;.
\end{multline*}

\bigskip

\noi By \eqref{comp:i-bra},  this is simply a consequence of
  \eqref{cocep5clas-b}, since
$\lam_{ijk} = i(g_{ijk})\;. \hspace{1.4cm}\Box$

\bigskip

 \noi We  now return to  our computation \eqref{cockap2a}:

\begin{eqnarray*}
\de^1_{\lambda_{ij}}(\de_{ij}) 
&=&  \de^1_{\lambda_{ij}}\de^0_{\lambda_{ij}}(B_i) 
- \de^1_{\lam_{ij}}\D^1_{m_i}(-\gamma_{ij})\\
&=&  \de^1_{\lambda_{ij}}\de^0_{\lambda_{ij}}(B_i)
 - d^1_{m_i} \de^1_{\lam_{ij}}(-\gamma_{ij})
\\ 
&=& [g_{ijk},\, B_i] - d^1_{m_i}(g_{ijk}\,d_{m_i}(g_{ijk}^{-1})) \\
&=& [g_{ijk},\,i_{B_i} -  dm_i]  \hspace{2.5cm} \text{by \eqref{dndn1a}}
\\
&=& [\nu_i, \, g_{ijk}]\,. 
\end{eqnarray*}
This finishes the  second proof of equation \eqref{cockap2} .
\begin{flushright}
$\Box$
\end{flushright}

 We now set 

\bee
\label{defom}
\om_i:= \D^2_{m_i}(B_i)\,.
\ee

\bigskip

\noi  Since the combinatorial definition of the twisted de Rham differential
 $\D^2$  (\cite{cdf} (3.3.1)) matches the global geometric  definition
 \eqref{cube} of the 3-curvature  $\Om$, this 3-curvature $\Om$ is
 locally described by the $G$-valued  3-forms $\om_i$.

\bigskip

It follows from the definitions \eqref{ifi} and \eqref{defom} of the
forms $\nu_i$ and $\om_i$, and from \eqref{dndn1},
that 
\begin{eqnarray*}
\D^3_{m_i} (\om_i) & = & \, \D^3_{m_i}\D^2_{m_i}(B_i)\\
& = &  [\D^1m_i,\, B_i]\\
& = & [\nu_i,\, B_i] +[B_i,\, B_i] 
\end{eqnarray*} 
so that  the local 3-curvature form $\om_i$ satisfies 
the higher Bianchi identity 
\bee
\label{relnufi}
\D^3_{m_i}(\om_i) = [\nu_i, B_i]\,.
\ee
A second  relation between the forms $\nu_i$ and 
$\om_i$  follows from their definitions and 
  the Bianchi identity for the
 $\mathrm{Aut}(G)$-valued  1-form   $m_i$:
\begin{eqnarray*}
i(\om_i) & = & \, \D^2_{m_i}i(B_i) \\
& = & \, \D^2_{m_i}(\D^1m_i - \nu_i)\\
& = & \,\D^2_{m_i} ( - \nu_i)\,,
\end{eqnarray*} in other words
\bee
\label{ificonj} 
\D^2_{m_i}\nu_i + i(\om_i) = 0\,.
\ee
This equation is  the local form of equation\eqref{eqom1}, just as
\eqref{relnufi} was the local form of \eqref{eqom2}.

\bigskip

We will now   show that  the  equation
\eqref{eq:bij1a}
for the 2-forms $B_i$, which we write here as
 \[
\label{eq:bij1b}
\de^0_{\lambda_{ij}}(B_i) = \D^1_{m_i}(-\ga_{ij}) + \de_{ij}\;,
\]
 induces the corresponding  gluing equation for
the local 3-forms $\om_i$. 
From the definition of $\lambda_{ij}(\om_j)$ and \eqref{funct:dnm}, it
follows that   
\begin{eqnarray*}
\lambda_{ij}(\om_j) & = & \lambda_{ij}(\D^2_{m_j}(B_j))\\
&=& \D^2_{{}^{\lambda_{ij}\,\ast\,}\!m_j}\lambda_{ij}(B_j)
\end{eqnarray*}
and by the gluing laws \eqref{cocep13clas} and \eqref{eq:bij1a} for
$m_i$ and $B_i$, this can be stated as
\begin{eqnarray*}
\lambda_{ij}(\om_j) &=& \D^2_{m_i - i(\ga_{ij})}(B_i + \de_{ij}
 + \D^1_{m_i}(-\,\ga_{ij}))\\
&=&
\D^2_{m_i}(B_i)  + \D^2_{m_i}(\de_{ij}) + \D^2_{m_i}\D^1_{m_i}(-\,\ga_{ij})
- [\ga_{ij},\, B_i + \de_{ij}+\D^1_{m_i}(-\,\ga_{ij})]\;.  
\end{eqnarray*}
By  \eqref{falsebianchi}, this last  equality can be rewritten as
\begin{eqnarray*}
\lambda_{ij}(\om_j) &=&\om_i  + \D^2_{m_i}(\de_{ij}) + [ \D^1m_i,\,\,-\,
\ga_{ij}] -  [\ga_{ij},\, B_i] -  [\ga_{ij},\,\de_{ij}]\\
&=& \om_i + \D^2_{m_i}(\de_{ij}) + [\ga_{ij},\, \D^1m_i - B_i]
 -  [\ga_{ij},\,\de_{ij}]
\end{eqnarray*}
and by \eqref{cockap01} this proves the gluing law for the 3-forms
   $\om_i$ \cite{dgg} (6.1.23):
\begin{eqnarray*}
\label{comoioj}
\lambda_{ij}(\om_j) =  \om_i + \D^2_{m_i}(\de_{ij}) +[\ga_{ij},\,
\nu_i] - [\ga_{ij},\, \de_{ij}]\;.
 \end{eqnarray*}
 By combining this  with the gluing law
\eqref{cockap1} for $\nu_i$, we see that \eqref{comoioj}  can finally be rewritten in the more
compact form 
\begin{equation}
\label{comoioj1}
\lambda_{ij}(\om_j) + [{}^{\lambda_{ij}\,}\!\nu_j, \, \ga_{ij}]
= \om_i + \D^2_{m_i}(\de_{ij})
\ee

\bigskip

\section{\v{C}ech-de Rham coboundaries}
 \label{section5}

\subsection{}
We saw  in section \ref{section1} how  a  change
 in the choice
trivializing data $(x_i,\,\phi_{ij})$  in a gerbe $\pc$ could be 
measured by a pair $(r_i,\, \theta_{ij})$ 
\eqref{def-mui},\eqref{defteij} inducing a coboundary relation  
\eqref{eq:cob3} between the corresponding cocycle pairs
$(\lambda_{ij},\, g_{ijk})$. We will now examine how the terms
$(m_i,\, \ga_{ij})$, $(\nu_i,\, \de_{ij})$ and $B_i$ introduced in
section \ref{section3} vary when the  arrows $\ga_i$ 
\eqref{def:gai} and $\de_i$ \eqref{def:dei} which determine them have
been   modified.

\bigskip

The difference between the  arrow  $\ga_i$ and 
an analogous  arrow
$\ga_i'$  is measured by a
1-form $e_i \in \LG \ot\, \Om^1_{U_i}$,
 defined by the commutativity of
the following diagram:

\begin{equation}
\label{def:ei}
\xymatrix{
\epsilon \prb x_i \ar[rr]^{\epsilon \prb \chi_i} \ar[d]_{\ga_i} &&
\epsilon \prb x'_i \ar[d]^{\ga'_i}\\
\pra  x_i \ar[r]_{\pra \chi_i} & \pra x'_i \ar[r]_{e_i}& \pra x'_i
}
\end{equation}
This conjugates to a commutative  diagram

\[
\xymatrix{
\prb G \ar[rr]^{\prb r_ i} \ar[d]_{m_i} && \prb G \ar[d]^{m'_i}\\
\pra G \ar[r]_{\pra r_i} & \pra G \ar[r]_{i(e_{i})} & \pra G
}
\]
so that 
\begin{eqnarray*}
m'_i &=& i(e_i) \, (\pra r_i) \, m_i \,\, (\prb r_i)^{-1}  \\
&=& i(e_i)\, [\pra r_i \,m_i \, {\pra r_i}^{-1}]\, [\pra r_i\, {\prb
  r_i}^{-1}]
\end{eqnarray*}
In classical terms, this
  is expressed
 as an equation
\begin{eqnarray}
m'_i &=& {}^{r_i\,} m_i + r_i \D r_i^{-1} + i(e_i) \label{eq:mi}    \\
    &=& {}^{r_i \ast\,}m_i + i(e_i)\;. \label{eq:mi1}
\end{eqnarray}
which compares the connections $m_i$ and $m'_i$ induced on the 
group $G$ by the arrows $\ga_i$ and $\ga'_i$.

\newpage

We now consider the following diagram in $\pc_{U_{ij}}$:

\bee
\label{comp:eij}
\xymatrix@C=40pt{
\pra x'_i \ar[d]_{e_i} & \pra x'_i \ar[l]_{ r_i(\ga_{ij})}
\ar[d]^{\pra \theta_{ij}}\\
\pra x'_i \ar[d]_{m'_i(\prb \theta_{ij})} & \pra x'_i
\ar[d]^{ \lambda'_{ij}(e_j)} \\
\pra x'_i & \ar[l]^{\gamma'_{ij}}\pra x'_i\;.
} 
\end{equation}

\begin{proposition}
\label{diag:gagaprime}
The diagram \eqref{comp:eij}  is commutative. 
\end{proposition}
\noi {\bf Proof:} Consider the diagram 
\begin{equation}
\label{diag:ei}
\xymatrix@R=15pt@C=28pt{
&&\epsilon \prb x_j \ar[ddll]_{\epsilon \prb \phi_{ij}}
\ar[rrrr]^{\ga_j}\ar@{-}[dd] 
&&&&\pra x_j \ar[ddd]_{\pra \chi_j} \ar[dl]^{\pra \phi_{ij}} \\
&&&&&\pra x_i \ar[dl]_{\gamma_{ij}}\ar[dd]^{\pra \chi_i}&\\
\epsilon \prb x_i \ar[ddddd]_{\epsilon \prb \chi_i }
\ar[rrrr]^(.6){\ga_{i}} &&\ar[dddd]^(.5){\epsilon \prb \chi_j}
&&\pra x_i
\ar[dd]_(.3){\pra \chi_i}  &&\\
&&&&&\pra x'_i\ar[d]^(.4){\pra \theta_{ij}}\ar[dl]_(.4){
\pra  r_i(\ga_{ij})}
&\pra x'_j\ar[dl]^{\pra \phi'_{ij}}
\ar[ddd]^{e_j}\\
&&&&\pra x'_i \ar[ddd]^(.4){e_i}
&\pra x'_i \ar[ddd]^(.3){\pra \lambda'_{ij}(e_j)}
 &\\
&&&&&&\\
&&\epsilon\pra x'_j \ar[ddll]_(.3){\epsilon \prb \phi'_{ij}}
|!{[dll];[drr]}\hole
 \ar@{-}[rr]^(.7){\gamma'_j}&& \ar@{-}[r]     
& \ar[r]&\pra x'_j\ar[dl]^{\pra \phi'_{ij}}\\
   \epsilon\prb x'_i \ar[rrrr]^{\ga'_i}
\ar[d]_{\epsilon \prb \theta_{ij}}
&&&&\pra x'_i \ar[d]_{m'_i(\prb \theta_{ij})} &\pra x'_i \ar[dl]^{\ga'_{ij}} &\\
\epsilon \prb x'_i  \ar[rrrr]_{\ga'_i} 
&&&&  \pra x'_i  &&
}\end{equation}
 The lower front square of the right-hand face of this cube
is just the square \eqref{comp:eij}.
  Since we know that all the other squares in this diagram commute,
 so does the square 
\eqref{comp:eij}. $\hspace{10.5cm}\qed$

\bigskip

\noi The  commutativity of \eqref{comp:eij} is equivalent to  the equation
\begin{equation}
\label{eq:eiteij}
m'_i(\prb \theta_{ij})\,e_i\,  r_i(\ga_{ij}) = \ga'_{ij}\, 
\lambda'_{ij}(e_j) \,  \pra  \theta_{ij}\;.
\end{equation}
This may be rewritten in classical notation as:
\bee
\label{eq:eiteij4a} 
(\ga'_{ij} - {}^{\theta_{ij}\,}\!r_i(\ga_{ij})) \,
+ \, (\lambda'_{ij}(e_j) - {}^{\theta_{ij}\,}\!e_i)
= \D_{m'_i}\,\th_{ij}\, \,\th_{ij}^{-1} \;.
\ee

\bigskip

We now choose a family of  arrows $\delta'_i: \kappa \pra x'_i \la
\pra x'_i$. The families $\de'_i$ and $\ga'_i$ determine  as in
\eqref{def:bi} a family  of
$\mathfrak{g}$-valued 2-form $B'_i$ above $U_i$. The latter  in turn
determines, together with the pair of form $(m'_i,\, \ga'_{ij})$
\eqref{eq:mi}, \eqref{eq:eiteij4a},
 a new pair of 2-forms $(\nu'_i,\, \de'_{ij})$  and a 3-form $\om'_i$
 satisfying the corresponding equations \eqref{cockap1},
 \eqref{ificonj},  \eqref{cockap2},  
 \eqref{relnufi} and  \eqref{comoioj1}. The families  $\de_i$ and $\de'_i$ are
 compared by the following analogue of diagram \eqref{def:ei}:
\begin{equation}
\label{def:ni}
\xymatrix{
\kappa \pra x_i \ar[rr]^{\kappa \pra \chi_i} \ar[d]_{\de_i} &&
\kappa \pra x'_i \ar[d]^{\de'_i}\\
\pra  x_i \ar[r]_{\pra \chi_i} & \pra x'_i \ar[r]_{n_i}& \pra x'_i\;\;.
}
\end{equation}

\noi We will now compare the 2-forms $B_i$ and $B'_i\,.$ 
We 
consider the   diagram

\begin{equation}
\label{bibi}
\xymatrix@C=0pt@R=4pt{
&&&\epsilon_{01} \epsilon_{12}(\prc
x_i)\ar[dddlll]_{\epsilon_{01}\epsilon_{12}(\prc \chi_i)} \ar[rrrrrr]^{K(\prc x_i)}
\ar@{-}[ddd]_(.7){\epsilon_{01}(\ga_i^{12})}&&&&&& 
\ka \epsilon_{02} (\prc x_i)  \ar[dddlll]_{\kappa\,
  \epsilon_{02}(\prc \chi_i)}\ar[dddd]_{\ka(\ga^{02}_i)}\\
&&&&&&&&&\\
&&&&&&&&&\\
\epsilon_{01} \epsilon_{12}(\prc x'_i)\ar[dddd]^(.4){\epsilon_{01}\,{\ga'}_i^{12}}
\ar[rrrrrr]_(.4){K(\prc x'_i)}
&&&\ar[d] &&&\ka \,\epsilon_{02} (\prc x'_i)\ar[dddd]_{\ka({\ga'}_i^{02})} &&&\\
&&&\epsilon_{01}\prb x_i \ar[ddll]_{\epsilon_{01}\prb \chi_i}
\ar[dddddd]^{\ga_i^{01}}
 &&&&&&\ka\,\pra x_i\ar[ddll]^{\ka
  \pra \chi_i}\ar[dddddd]_{\de_i}
\\
&&&&&&&&&\\
&\epsilon_{01}\prb x'_i\ar[dl]
 \ar[dddddd]^{{\ga'}_i^{01}} &&&&  & 
&   \ka\pra x'_i\ar[dddddd]_{\de'_i}
\ar[dl]^{\ka (e_{i}^{02})} &\\
\epsilon_{01}\,\pra \,x'_i \ar[dddddd]^{{\ga'}_i^{01}}
&&&&&&\ka \pra x'_i \ar[dddddd]_(.3){\de'_i} & &&\\
&&&&&&&&&\\
&&&&&&&&&\\
&&& \pra x_i  \ar@{<-}[rrr]^(.6){B_i}\ar[dl]_(.5){\pra \chi_i} &&&
\ar@{-}[r]&
\ar@{-}[rr]&&\pra x_i\ar[dl]^{\pra \chi_i}\\
&&\pra x'_i \ar@{<-}[rrrr]^(.6){r_i(B_i)} \ar[dl]^{e_{01}^i} &&&&
\ar@{-}[r]&\ar@{-}[r]
&\pra x'_i \ar[dl]^{n_i} &\\
&\pra x'_i \ar[dl]
&&&&&&\pra x'_i \ar[dl]^{\nu'_i(e_i^{02})}  &&\\
\pra x'_i \ar@{<-}[rrrrrr]_{B'_i}
&&&&&&\pra x'_i&&&
}\end{equation}

\bigskip

\noi in which the upper and lower unlabelled arrows are respectively
 $\epsilon_{01}(\prb e_i^{12})$ and $m'_i{}^{01}(e_i^{12})$. 

\newpage
\noi  The front  square (or rather hexagon) of the bottom face
\[
\xymatrix@C=80pt{
\pra x'_i \ar@{<-}[r]^{r_i(B_i)} \ar[d]_{e^i_{01}} & \pra x'_i 
\ar[d]^{\pra n_i} \\
\pra x'_i \ar[d]_{m'_i{}^{01}(e_i^{12})} & \pra x'_i
  \ar[d]^{\nu'_i(e^{02}_i)}\\
\pra x'_i \ar@{<-}[r]_{B'_i} &\pra x'_i 
}
\]
is commutative, since  all other squares in diagram
\eqref{bibi}
 are. Equivalently, since the action of the $\mathrm{Aut}(G)$-valued
 2-form  $\nu'_i$ on $e_i^{02}$ is
 trivial, this proves that  the equation 
\bee
\label{eq:bbprime}
 B'_i = r_i(B_i)  - \D^1_{m'_i}(-e_i) - n_i\,.
\ee
is satisfied. In particular for given $B_i$ and $e_i$, the 2-forms
 $B'_i$ and $n_i$ actually determine each other. 

\vspace{1cm}

By conjugation, diagram \eqref{def:ni} induces a commutative diagram
\[
\xymatrix{
\pra G \ar[rr]^{\pra r_ i} \ar[d]_{\nu_i} && \pra G \ar[d]^{\nu'_i}\\
\pra G \ar[r]_{\pra r_i} & \pra G \ar[r]_{i_{n_{i}}} & \pra G
}
\]
equivalent to the equation
\[ i(n_i)\; \pra r_i \; \nu_i = \nu'_i \; \pra r_i\;.\]
In classical terms, this is the simpler  analogue
\bee
\label{eq:ni1}
\nu'_i =   {}^{r_i\,}\!\nu_i +  i(n_i)
\end{equation}
 for $\nu_i$ 
of  the equation \eqref{eq:mi} for $m_i$.

\bigskip

We will now show that this coboundary equation for $\nu_i$
 can be derived from the definition \eqref{ifi} of
$\nu_i$, and the coboundary equations \eqref{eq:mi}
 and \eqref{eq:bbprime} for
$m_i$ and $B_i$:

\begin{eqnarray*}
\nu'_i &=&  \D^1m'_i - i(B'_i)\\
&=& \D^1({}^{r_i \ast\,}\!m_i +i(e_i)) - i(r_i(B_i)+ n_i
+\D^1_{m'_i}(-e_i)) \\ &=&
{}^{r_i\,}\!\D^1m_i + i(\D^1e_i) +[{}^{r_i\ast\,}\!m_i,\,i(e_i)] -
i(r_i(B_i)) + i(\D^{1}_{m'_i}(-e_i)) +i(n_i)\\ &=&
{}^{r_i\,}\!(\D^1m_i - i(B_i) ) + i(n_i)  + i(\D^{1}_{m'_i}(-e_i)
+ \D^1e_i + [{}^{r_i \ast\,}\!m_i,\,e_i])
\end{eqnarray*}
In order to prove \eqref{eq:ni1}, it now suffices to verify that the 3
terms in 
 the last summand of the final  equation cancel each
other out:

\begin{eqnarray*}
\D^1_{m'_i} (-e_i) + \D^1(e_i) +   [{}^{r_i \ast\,}\!m_i,\,e_i] 
&=& \D^1(-e_i) - [m'_i,\,e_i] +\D^1e_i + [{}^{r_i
  \ast\,}\!m_i,\,e_i] \\
&=& \D^1(-e_i) + \D^1e_i -[e_i,\,e_i] \\
&=& 0 \;.  \hspace{6cm}\Box
\end{eqnarray*}

The other equation satisfied by the forms $n_i$ is the counterpart of equation
\eqref{eq:eiteij}. It is obtained by considering the following  diagram,
analogous to \eqref{diag:ei}:

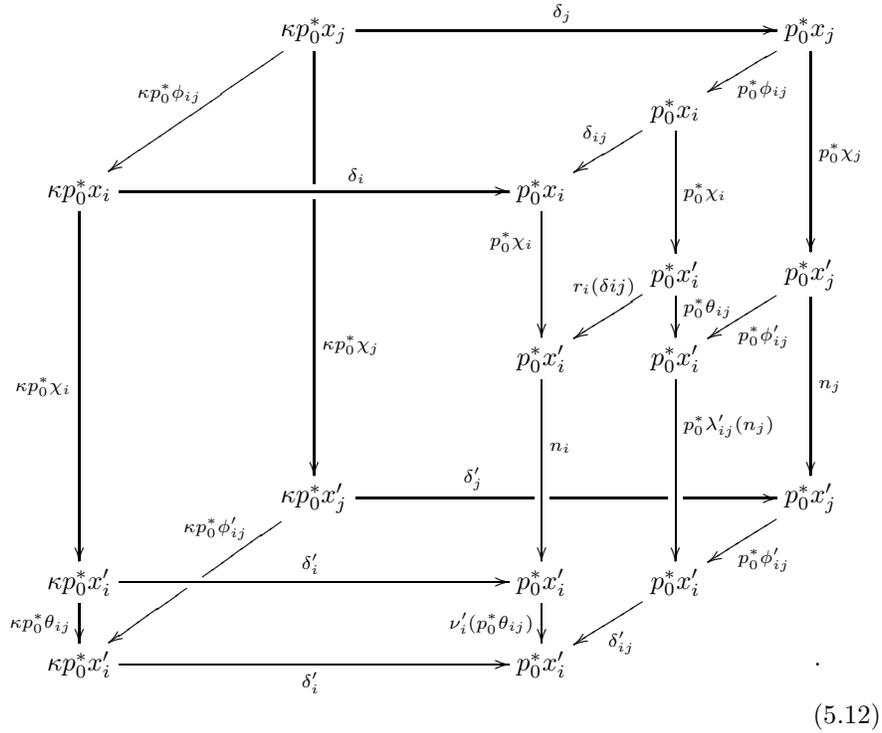
\bee
\label{diag:ni}
\xymatrix@R=15pt@C=26pt{
&&\kappa \pra x_j \ar[ddll]_{\kappa \pra \phi_{ij}}
\ar[rrrr]^{\de_j}\ar@{-}[dd] 
&&&&\pra x_j \ar[ddd]^{\pra \chi_j} \ar[dl]^{\pra \phi_{ij}} \\
&&&&&\pra x_i \ar[dl]_{\de_{ij}}\ar[dd]^{\pra \chi_i}&\\
\kappa \pra x_i \ar[ddddd]_{\kappa \pra \chi_i }
\ar[rrrr]^(.6){\de_{i}} &&\ar[dddd]^(.5){\kappa \pra \chi_j}
&&\pra x_i
\ar[dd]_(.3){\pra \chi_i}  &&\\
&&&&&\pra x'_i\ar[d]^(.4){\pra \theta_{ij}}\ar[dl]_(.35){r_i(\de{ij})}
&\pra x'_j\ar[dl]^{\pra \phi'_{ij}}
\ar[ddd]^{n_j}\\
&&&&\pra x'_i  \ar[ddd]^(.4){n_i}
&\pra x'_i \ar[ddd]^(.3){\pra \lambda'_{ij}(n_j)}
 &\\
&&&&&&\\
&&\kappa\pra x'_j \ar[ddll]_(.3){\kappa \pra \phi'_{ij}} |!{[lld];[drr]}\hole
 \ar@{-}[rr]^(.7){\de'_j}&& \ar@{-}[r]     
& \ar[r]&\pra x'_j\ar[dl]^{\pra \phi'_{ij}}\\
   \kappa\pra x'_i \ar[rrrr]^{\de'_i}
\ar[d]_{\kappa \pra \theta_{ij}}
&&&&\pra x'_i \ar[d]_{\nu'_i(\pra \theta_{ij})} &\pra x'_i \ar[dl]^{\de'_{ij}} &\\
\kappa \pra x'_i  \ar[rrrr]_{\de'_i} 
&&&&  \pra x'_i  &&\;\;.
}\end{equation}

The lower front square on the  right-hand face
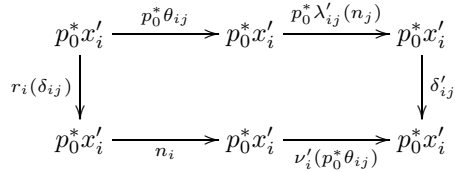
\[\xymatrix@C=40pt{
\pra x'_i \ar[r]^{\pra \theta_{ij}} \ar[d]_{r_i(\de_{ij})} & \pra x'_i 
\ar[r]^{\pra \lambda'_{ij}(n_j)} & \pra x'_i \ar[d]^{\de'_{ij}}\\
\pra x'_i \ar[r]_{n_i} & \pra x'_i \ar[r]_{\nu'_i(\pra \theta_{ij})}
 & \pra x'_i
}
\]
of diagram \eqref{diag:ni} is commutative, since all other squares
in this diagram are.

\pagebreak

\noi  This proves that  equation 
\[
\label{eq:ninj}
\xymatrix{
\nu'_i(\pra\, \theta_{ij}) \,\, n_i \,\, r_i(\de_{ij}) = 
\de'_{ij}\,\, \pra \lambda'_{ij}(n_j) \,\,\pra \theta_{ij} 
}
\]
in $\LG \ot \Om^2_{U_i/S}$ is satisfied. Regrouping the various terms
in this equation
as we did above for  equation \eqref{eq:eiteij},
 we find that it is equivalent, in additive notation, to

\[
\label{eq:niri}
(\de'_{ij} - r_i(\de_{ij})) + ( \lambda'_{ij}(n_j) -
{}^{\th_{ij}\,}\!n_i) = [\nu'_i,\, \th_{ij}]\;,
\]
an equation for 2-forms very  similar to equation
\eqref{eq:eiteij4a}
for 1-forms.

\bigskip

We will now examine  the effect of  the chosen transfomations
\bee
\label{eq:eiteij4b}
 (\lam_{ij},\, g_{ijk},\, m_i,\, \ga_{ij}) \qquad \la \qquad 
(\lam'_{ij} ,\,g'_{ijk},\, m'_i,\, \ga'_{ij}) \ee
 and  $B_i \la B_i'\ $ \eqref{eq:bbprime} on the 
3-curvature 3-forms $\om_i$ \eqref{defom}. For this, it will be convenient to
set 

\[
\label{eq:bar}
\bar{e}_i:= r_i^{-1}(e_i) \qquad \text{and} \qquad
 \bar{n}_i := r_i^{-1}(n_i)\;.
\]
It follows from \eqref{addm},  \eqref{compbra}, and  the
transformation formula  \eqref{eq:mi1} that
\bee
\label{eq:transdm}
\D^n_{m'_i}(r_i(\eta)) = r_i (\D^n_{m_i}(\eta) + [\bar{e}_i,\,
\eta])
\ee
for any $G$-valued $n$-form $\eta$  with $n >1$.
 In particular
\begin{eqnarray*}
\D^1_{m'_i}(-e_i) &=& \D^1_{\,\,{}^{r_i \ast}m_i}(-e_i) -[e_i,\, e_i]\\
&=& r_i(\D^1_{m_i}(-\bar{e}_i) -[\bar{e}_i,\, \bar{e}_i])
\end{eqnarray*}
 so that  \eqref{eq:bbprime} may be expressed as
\[
B'_i = r_i(B_i - \D^1_{m_i}(- \bar{e}_i )
+ [\bar{e}_i,\, \bar{e}_i]  - \bar{n}_i )\;.
\]
Applying once more the formula  \eqref{eq:transdm}, we   find that  
\begin{eqnarray}
\label{coboun-om}
\hspace{1cm}\om'_i & = &
\D^2_{m'_i}(B'_i) \notag \\
&= &\D^2_{m'_i}(    r_i(B_i - \D^1_{m_i}(- \bar{e}_i )
+ [\bar{e}_i,\, \bar{e}_i]  - \bar{n}_i )    ) \notag \\
&=&
 r_i (\D^2_{m_i}(B_i - \D^1_{m_i}(- \bar{e}_i )
+ [\bar{e}_i,\, \bar{e}_i]  - \bar{n}_i )) +
   \notag \\
&&   \hspace{2.2cm} +\, [\bar{e}_i,\, B_i -  \D^1_{m_i}(- \bar{e}_i )
+\, [\bar{e}_i,\, \bar{e}_i]  - \bar{n}_i]\;.
\end{eqnarray}
We now make use of   \eqref{falsebianchi} in order to compute 
the value of the  expression $ \D^2_{m_i}\D^1_{m_i}(-\bar{e}_i)$ which
arises  when we
expand the first summand of the last equation  \eqref{coboun-om}:
\begin{eqnarray*}
 \D^2_{m_i}\D^1_{m_i}(-\,\bar{e}_i) &=&[\D^1m_i\,,\,-\,\bar{e}_i] +
 [\D^1(-\,\bar{e}_i),\, -\, \bar{e}_i] +
 [[m_i,\,-\,\bar{e}_i],\,-\,\bar{e}_i]\\
&=&- \, [\D^1m_i\,,\, \bar{e}_i] +[\D^1\bar{e}_i,\,\bar{e}_i] +[[m_i,\,
\bar{e}_i],\, \bar{e}_i]\;.
\end{eqnarray*}
  Inserting this expression into \eqref{coboun-om}, we find
that
\begin{multline}
\label{omom1}
\om'_i = r_i(\om_i + [\D^1m_i,\, \bar{e}_i] - [\D^1\bar{e}_i,\,
\bar{e}_i] -[[m_i,\, \bar{e}_i],\, \bar{e}_i] - \D^2_{m_i}(\bar{n}_i)\; +
\\ 
+ \D^2_{m_i}[\bar{e}_i,\, \bar{e}_i] + [\bar{e}_i,\, B_i] -
[\bar{e}_i,\, \D^1_{m_i}(-\bar{e}_i)]
 -[\bar{e}_i,\, \bar{n}_i]  )\;.
\end{multline}
The four terms 
\[ - [\D^1\bar{e}_i,\,
\bar{e}_i]  -[[m_i,\, \bar{e}_i],\, \bar{e}_i]+
 \D^2_{m_i}[\bar{e}_i,\, \bar{e}_i] -[\bar{e}_i,\,
 \D^1_{m_i}(-\bar{e}_i)]
\]
cancel each other out, so that we are left in \eqref{omom1}
with 
\begin{eqnarray}
\label{omom2}
\om'_i&=& r_i(\om_i+ [\D^1m_i,\, \bar{e}_i] - \D^2_{m_i}(\bar{n}_i)
+  [\bar{e}_i,\, B_i] -[\bar{e}_i,\, \bar{n}_i]) \notag\\
&=& r_i(\om_i + [\D^1m_i  -i(B_i)   ,\,\bar{e}_i] + [\bar{n_i},\,
\bar{e}_i] - \D^2_{m_i}(\bar{n}_i))\notag\\
&=&  r_i(\om_i) + r_i([\nu_i  ,\,\bar{e}_i]) + 
r_i([ \bar{n_i}  ,\,\bar{e}_i])
 - r_i( \D^2_{m_i}(\bar{n}_i))\notag\\
&=&  r_i(\om_i) + [{}^{r_i\,}\! \nu_i,\, e_i] + [n_i,\, e_i] -
\D^2_{\,{}^{r_i\, \ast\,}\!m_i}(n_i)
\end{eqnarray} 
where in the last line  we made use of the functoriality property
\eqref{compbra} of the bracket operation. Amalgamating the last two
summands, we may finally write 
 the coboundary transformation for the 3-curvature form
$\om_i$ in the compact form
\begin{eqnarray*}
\label{coboun-om1}
\om'_i
&=&  r_i(\om_i) +  [{}^{r_i\,}\! \nu_i,\, e_i] - \D^2_{m'_i}(n_i)
\;.
\end{eqnarray*}
If instead we amalgamate the second and third term in \eqref{omom2},
we find the equivalent formulation

\begin{eqnarray}
\label{coboun-om1a}
\om'_i
&=&  r_i(\om_i) +[\nu'_i,\,e_i] -\D^2_{{}^{r_i\, \ast\,}\!m_i}(n_i)
\;.
\end{eqnarray}

\bigskip

\begin{lbremark}
\label{remark}{\rm {\bf(Comparison with \cite{dgg}):}

\bigskip

The coboundary equation \eqref{coboun-om1a} is compatible with equation (6.2.19) of \cite{dgg}, but neither is a
special case of the other. Here   we
allowed  both the trivializing data $(x_i,\, \phi_{ij})$ for the gerbe
and the  expressions $(\ga_i,\,\de_i,\, B_i)$ for the curving data to vary, whereas in the
coboundary equations of \cite{dgg} the gerbe data  $(x_i,\, \phi_{ij})$ was
fixed and only  the  $(\ga_i,\,\de_i,\, B_i)$  varied. This restriction
amounted  to setting $(r_i,\, \theta_{ij}) = (1,\, 1)$ in our
equation \eqref{eq:eiteij4a}. On the other hand,  a
 notion of equivalence between cocycles was introduced in
 \cite{dgg}
 which was more extensive than
 the one considered here. In order for these  to be comparable, one
 must 
 suppose that  the arrow $h$
in diagram (4.2.1) of \cite{dgg} is the
identity map, i.e.  that    the pair of
differential forms $(\pi_i,\, \eta_{ij})$ associated to $h$ in  {\it loc. cit} \S$\, 6.2$
is  trivial. This is a reasonable assumption, since a non-trivial arrow  $h$
could    really be termed a gauge transformation, rather than
a coboundary term. 
With this additional condition, the last two summands in equation
(6.2.19)
of \cite{dgg}  vanish, so that  this equation reduces to
\begin{equation}
\label{eq:simp}
\om'_i = \om_i +\delta^2_{m_i}(\alpha_i) -[\nu'_i,\, E_i]\;.
\end{equation}
This simplified equation  is  compatible with our equation
  \eqref{coboun-om1a} with $r_i= 1$,  under the correspondence
 $e_i := - E_i$ and
$n_i := -\, \alpha_i$. 
 }
\end{lbremark}


\end{document}